\documentclass[aap,MSNbibl,seceqn,dvips]{arximspdf}
\usepackage{graphicx}
\doi{10.1214/13-AAP979} %kopijuoti is PTS
\volume{24}
\issue{6}
\pubyear{2014}
\firstpage{2340}
\lastpage{2370}

\makeatletter
\newtheorem{prop}{Proposition}
\newtheorem{teo}{Theorem}
\newproclaim{rem}{Remark}
\newproclaim{exa}{Example}
\makeatother

\begin{document}
\begin{frontmatter}

\title{Quickest detection of a hidden target and~extremal surfaces}
\runtitle{Quickest detection of a hidden target}

\begin{aug}
\author{\fnms{Goran} \snm{Peskir}\corref{}\ead[label=e1]{goran@maths.man.ac.uk}}
\runauthor{G. Peskir}
\affiliation{University of Manchester}
\address{School of Mathematics\\
University of Manchester\\
Oxford Road\\
Manchester M13 9PL\\
United Kingdom\\
\printead{e1}} %adresu isvedimo komanda gale!
\end{aug}

% HISTORY:
\received{\smonth{1} \syear{2013}}

% ABSTRACT
%
\begin{abstract}
Let $Z=(Z_t)_{t \ge0}$ be a regular diffusion process started at
$0$, let $\ell$ be an independent random variable with a strictly
increasing and continuous distribution function $F$, and let
$\tau_\ell= \inf\{ t \ge0 \vert Z_t = \ell\}$ be the
first entry time of $Z$ at the level $\ell$. We show that the
quickest detection problem
\[
\inf_\tau \bigl[ \mathsf P(\tau<
\tau_\ell) + c \mathsf E (\tau - \tau_\ell)^+ \bigr]
\]
is equivalent to the (three-dimensional) optimal stopping problem
\[
\sup_\tau\mathsf E \biggl[ R_\tau- \int
_0^\tau c(R_t) \,dt \biggr],
\]
where $R = S - I$ is the range process of $X=2F(Z)-1$ (i.e., the
difference between the running maximum and the running minimum of
$X$ ) and $c(r) = c r$ with $c>0$. Solving the latter problem
we find that the following stopping time is optimal:
\[
\tau_* = \inf \bigl\{ t \ge0 \vert f_*(I_t,S_t)
\le X_t \le g_*(I_t,S_t) \bigr\},
\]
where the surfaces $f_*$ and $g_*$ can be characterised as extremal
solutions to a couple of first-order nonlinear PDEs expressed in
terms of the infinitesimal characteristics of $X$ and $c$. This is
done by extending the arguments associated with the maximality
principle [\textit{Ann. Probab.} \textbf{26} (1998) 1614--1640] to the three-dimensional setting of the
present problem and disclosing the general structure of the solution
that is valid in all particular cases. The key arguments developed
in the proof should be applicable in similar multi-dimensional
settings.
\end{abstract}

% KEYWORDS

\end{frontmatter}

%s1 #&#
\section{Introduction}\label{sec1}
%%%%%%%%%%%%%%%%%%%%%%

Imagine that you are observing a sample path $t \mapsto Z_t$ of the
continuous process $Z$ started at $0$ and that you wish to detect
when this sample\vadjust{\goodbreak} path reaches a level $\ell$ that is not directly
observable. Situations of this type occur naturally in many applied
problems, and there is a whole range of hypotheses that can be
introduced to study various particular aspects of the problem.
Assuming that $Z$ and $\ell$ are independent, and denoting by
$\tau_\ell$ the first entry time of $Z$ at $\ell$, it was shown
recently (see~\cite{Pe-5}) that the median/quantile rule minimises
not only the spatial expectation $\mathsf E [(\ell - X_\tau)^+ + c
(X_\tau - \ell)^+]$ (dating back to R. J. Boscovich 1711--1787)
but also the \emph{temporal} expectation $\mathsf E [(\tau_\ell -
\tau)^+ + c (\tau - \tau_\ell)^+ ]$ over all stopping times
$\tau$ of $Z$ where $c$ is a positive constant. Motivated by this
development, and seeking for further insights and connections, in
this paper we study the ``mixed'' variational problem
%
%e1.1 #&#
\begin{equation}
\label{11} \inf_\tau \bigl[ \mathsf P(\tau<
\tau_\ell) + c \mathsf E (\tau - \tau_\ell)^+ \bigr],
\end{equation}
which appears in the classic formulation of quickest detection due
to Shiryaev (see \cite{Sh-1,Sh-2} and \cite{PS}, Sections~22 and 24
and the references therein). The key difference between (\ref{11})
and the classic formulation is that the unobservable time
$\tau_\ell$ in~(\ref{11}) is obtained through the uncertainty in
the space domain (as the first entry time of $Z$ at the unknown
level $\ell$), while the unobservable time in the classic
formulation is obtained through the uncertainty in the time domain
(as the unknown level itself). Unlike the classic formulation,
however, we do not assume that the probabilistic characteristics of
$Z$ change following $\tau_\ell$ so that there is no learning about
the position of $\ell$ through the observation of $Z$ (quickest
detection problems of this kind require a different treatment and
will be studied elsewhere). Likewise, since the underlying loss
processes $t \mapsto1(t < \tau_\ell)$ and $t \mapsto1(t -
\tau_\ell)^+$ are not adapted to the natural filtration generated by
$Z$ (or its usual augmentation), we see that problem (\ref{11})
belongs to the class of ``optimal prediction'' problem (within optimal
stopping). Similar optimal prediction problems have been studied in
recent years by many authors (see, e.g., \cite{BDP,Co,DP-1,DP-2,DP-3,DPS,EE,ET,GHP,GPS,NS,Ped-2,Sh-3,Sh-4,Ur}). It may be noted in this
context that the nonadapted factor $\tau_\ell$ in the optimal
prediction problem (\ref{11}) is not revealed at the ``end'' of time
(i.e., it is not measurable with respect to the $ \sigma$-algebra
generated by the process $Z$).

While the median/quantile rule was derived in \cite{Pe-5} for
general (continuous) processes, a closer analysis of the mixed
variational problem (\ref{11}) reveals that this generality can
hardly be maintained. For this reason we restrict our attention to a
smaller class of processes and assume that $Z=(Z_t)_{t \ge0}$ is a
one-dimensional diffusion starting at $0$ and solving
%
%e1.2 #&#
\begin{equation}
\label{12}
dZ_t = a(Z_t) \,dt + b(Z_t)\,dB_t,
\end{equation}
where $a$ and $b>0$ are continuous functions, and $B = (B_t)_{t \ge
0}$ is a standard Brownian motion. To gain tractability we also
assume that the distribution function $F$ of $\ell$ is strictly
increasing and twice continuously differentiable. In the first step
we show that problem (\ref{11}) is equivalent to the optimal
stopping problem
%
%e1.3 #&#
\begin{equation}
\label{13} \sup_\tau\mathsf E \biggl[ R_\tau-
\int_0^\tau c(R_t) \,dt \biggr],
\end{equation}
where $R = S - I$ is the range process of $X = 2 F(Z) - 1$ (i.e., the
difference between the running maximum and the running minimum
of $X$) and $c(r) = c r$. This problem is of independent
interest and the appearance of the range process is novel in this
context revealing also that the problem is fully three-dimensional.
Two-dimensional versions of a related problem (when $I \equiv0$ and
$c$ constant) were initially studied and solved in important special
cases of diffusion processes in \cite{DS,DSS} and \cite{Ja}.
The general solution to problems of this kind was derived in the
form of the maximality principle in \cite{Pe-1}; see also Section~13
and Chapter~V in~\cite{PS} and the other references therein. In
these two-dimensional problems $c$ was a function of $X_t$ instead.
More recent contributions and studies of related problems include
\cite{CHO,Ga-1,Ga-2,GZ,Ho,Ob-1,Ob-2,Ped-1}; see also \cite{BDR,BI,HS}
and \cite{Pe-2} for related results in optimal control theory.
Close three-dimensional relatives of the problem (\ref{13}) also
appear in the recent papers \cite{DGM} and \cite{Zi} where the
problems were effectively solved by guessing and finding the optimal
stopping boundary in a closed form. These optimal stopping
boundaries are still curves in the state space.

In this paper we show how problem (\ref{13}) can be solved when
(i) no closed-form solution for the candidate stopping boundary is
available, and (ii) the optimal stopping boundaries are no longer
curves in the state space. This is done by extending the arguments
associated with the maximality principle \cite{Pe-1} to the
three-dimensional setting of the problem (\ref{13}) and disclosing
the general structure of the solution that is valid in all
particular cases. In this way we find that that the optimal stopping
boundary consists of two surfaces which can be characterised as
extremal solutions to a couple of first-order nonlinear PDEs. More
precisely, replacing $c(r)$ in problem (\ref{13}) above with a
more general function $c(i,x,s)$ specified below, we show that the
following stopping time is optimal:
%
%e1.4 #&#
\begin{equation}
\label{14} \tau_* = \inf \bigl\{ t \ge0 \vert f_*(I_t,S_t)
\le X_t \le g_*(I_t,S_t) \bigr\},
\end{equation}
where the surfaces $f_*$ and $g_*$ can be characterised as the
minimal and maximal solutions to
%
%e1.5 #&#
%e1.6 #&#
\begin{eqnarray}
\frac{\partial f}{\partial i}(i,s) &=& \frac{(\sigma^2 /2)(f(i,s))
L'(f(i,s))}{c(i,f(i,s),s) [L(f(i,s)) - L(i)]}
\nonumber\\[-8pt]\label{15} \\[-8pt]
&&{}\times \biggl[ 1 - \int_i^{f(i,s)} \frac{\partial c}{\partial i}(i,y,s)
\frac{L(y) - L(i)}{(\sigma^2 /2)(y) L'(y)} \,dy \biggr],\nonumber
\\
\frac{\partial g}{\partial s}(i,s) &=& \frac{(\sigma^2
/2)(g(i,s)) L'(g(i,s))}{c(i,g(i,s),s) [L(s) - L(g(i,s))]}
\nonumber\\[-8pt]\label{16} \\[-8pt]
&&{}\times  \biggl[ 1 + \int
_{g(i,s)}^s \frac{\partial c}{\partial s}(i, y,s)
\frac{L(s) - L(y)}{ (\sigma^2 /2)(y) L'(y)} \,dy \biggr]\nonumber
\end{eqnarray}
staying strictly above/below the lower/upper diagonal in the state
space, respectively (Theorem~\ref{teo1}). In these equations $\sigma$ is the
diffusion coefficient and $L$ is the scale function of $X$. They can
be expressed explicitly in terms of $a$, $b$ and~$F$. Recalling that
problems (\ref{11}) and (\ref{13}) are equivalent, we see that
this also yields the solution to the initial problem (\ref{11}). A
plain comparison with the median/quantile rule from \cite{Pe-5}
shows that the structure of problem (\ref{11}) is inherently
more complicated and the optimal stopping time $\tau_*$ may be
viewed as a nonlinear median/quantile rule. The optimal surfaces
$f_*$ and $g_*$ combined with the excursions of $X$ away from $I$
and $S$ exhibit interesting dynamics (not present in the
two-dimensional setting) which we describe in fuller detail as we
progress below. This dynamics may be combined with Lagrange
multipliers to tackle the constrained variant of the problem
(\ref{11}) where the probability error of early stopping is bounded
from above (we do not pursue this in the present paper). It is also
easily seen that swapping the order of $\tau$ and $\tau_\ell$ in
(\ref{11}) leads to optimal stopping at the diagonal and thus
corresponds to the linear median/quantile rule. The key arguments
developed in the proof rely heavily upon the extremal properties of
the optimal surfaces and should be applicable in similar
multi-dimensional settings.

%s2 #&#
\section{Quickest detection of a hidden target}\label{sec2}
%%%%%%%%%%%%%%%%%%%%%%%%%%%%%%%%%%%%%%%%%%%%%%%

In this section we will first formulate the quickest detection of
a hidden target problem and then show that this problem is
equivalent to an optimal stopping problem for the range process. The
latter problem will be studied in the next section.

Let $Z=(Z_t)_{t \ge0}$ be a one-dimensional diffusion process
starting at $0$ and solving
%
%e2.1 #&#
\begin{equation}
\label{21} dZ_t = a(Z_t) \,dt + b(Z_t)
\,dB_t,
\end{equation}
where $a$ and $b>0$ are continuous functions, and $B = (B_t)_{t \ge
0}$ is a standard Brownian motion. To meet a sufficient condition
used in the proof of Theorem~\ref{teo1} below we will also assume that $b^2$
is (locally) Lipschitz. Let $\ell$ be an independent random variable
with values in $\mathbb{R}$, and let
%
%e2.2 #&#
\begin{equation}
\label{22} \tau_\ell= \inf\{ t \ge0 \vert Z_t = \ell\}
\end{equation}
be the first entry time of $Z$ at the level $\ell$. We consider the
quickest detection problem
%
%e2.3 #&#
\begin{equation}
\label{23} V_1 = \inf_\tau \bigl[ \mathsf P(
\tau< \tau_\ell) + c \mathsf E (\tau - \tau_\ell)^+ \bigr],
\end{equation}
where the infimum is taken over all stopping times $\tau$ of $Z$
[i.e., with respect to the natural filtration $(\mathcal{F}_t^Z)_{t
\ge
0}$ generated by $Z$], and $c>0$ is a given and fixed constant (note
that whenever we say a stopping time throughout we always mean a
finite valued stopping time). Note that $\mathsf P(\tau< \tau_\ell)$
represents the probability of early stopping and $\mathsf E (\tau -
\tau_\ell)^+$ represents the expectation of late stopping when a
stopping time $\tau$ of $Z$ is being applied. Our task therefore is
to minimise the weighted sum of both errors over all stopping times
$\tau$ of $Z$. Note that $\ell$ and $\tau_\ell$ are not observable.
Set
%
%e2.4 #&#
\begin{equation}
\label{24} I_t^Z = \inf_{0 \le s \le t}
Z_s\quad\mbox{and}\quad S_t^Z = \sup
_{0
\le s \le t} Z_s
\end{equation}
for $t \ge0$, and let $F$ denote the distribution function of
$\ell$.

%pr1 #&#
\begin{prop}\label{prop1}
Problem (\ref{23}) is equivalent to the optimal stopping problem
%
%e2.5 #&#
\begin{equation}
\label{25} V_2 = \sup_\tau\mathsf E \biggl[ F
\bigl(S_\tau^Z \bigr) - F \bigl(I_\tau^Z-
\bigr) - c \int_0^\tau \bigl[ F
\bigl(S_t^Z \bigr) - F \bigl(I_t^Z-
\bigr) \bigr] \,dt \biggr],
\end{equation}
where the infimum is taken over all stopping times $\tau$ of $Z$.
\end{prop}

\begin{pf}
Let a stopping time $\tau$ of $Z$ be given and
fixed. First, using that $\ell$~and~$Z$ are independent, we find
that
%
%e2.6 #&#
\begin{eqnarray}\label{26}
\mathsf P(\tau< \tau_\ell) &=& 1 - \mathsf P(\tau\ge
\tau_\ell)\nonumber
\\
& =& 1 - \mathsf P(\tau \ge\tau_\ell, \ell> 0) -
\mathsf P(\tau\ge\tau_\ell, \ell \le0)\nonumber
\\
&=& 1 - \mathsf P \bigl(S_\tau^Z \ge\ell> 0
\bigr) - \mathsf P \bigl(I_\tau^Z \le\ell\le0 \bigr)
\\
&=& 1 - \mathsf E F \bigl(S_\tau^Z \bigr) + \mathsf E F
\bigl(I_\tau^Z- \bigr)\nonumber
\\
&=& 1 - \mathsf E \bigl[ F \bigl(S_\tau^Z \bigr)
- F \bigl(I_\tau^Z- \bigr) \bigr].\nonumber
\end{eqnarray}
Second, using a well-known argument (see, e.g., \cite{PS}, page 450)
it follows that
%
%e2.7 #&#
\begin{eqnarray}
\label{27} \mathsf E (\tau - \tau_\ell)^+ &=& \mathsf E \int
_0^\tau1(\tau _\ell\le t) \,dt = \mathsf E
\int_0^\infty1(\tau_\ell\le t) 1(t < \tau)
\,dt\nonumber
\\
&=& \int_0^\infty\mathsf E \bigl[ \mathsf
E \bigl( 1(\tau_\ell\le t) 1(t < \tau) \vert\mathcal{F}_t^Z
\bigr) \bigr] \,dt
\nonumber\\[-8pt]\\[-8pt]
&=& \int_0^\infty \mathsf E \bigl[ 1(t <
\tau) \mathsf E \bigl( 1(\tau_\ell\le t) \vert\mathcal
{F}_t^Z \bigr) \bigr] \,dt\nonumber
\\
&=& \mathsf E \int_0^\tau\mathsf P
\bigl(\tau _\ell\le t \vert\mathcal{F}_t^Z
\bigr) \,dt.\nonumber
\end{eqnarray}
Moreover, since $\ell$ and $Z$ are independent, we see that
%
%e2.8 #&#
\begin{eqnarray}
\label{28} \mathsf P \bigl(\tau_\ell\le t \vert
\mathcal{F}_t^Z \bigr) &=& \mathsf P \bigl(
\tau_\ell \le t, \ell> 0 \vert\mathcal{F}_t^Z
\bigr) + \mathsf P \bigl(\tau _\ell\le t, \ell\le0 \vert
\mathcal{F}_t^Z \bigr)\nonumber
\\
&=& \mathsf P \bigl(S_t^Z \ge \ell> 0 \vert
\mathcal{F}_t^Z \bigr) + \mathsf P \bigl(I_t^Z
\le\ell\le0 \vert \mathcal{F}_t^Z \bigr)
\\
&=& F \bigl(S_t^Z \bigr) - F
\bigl(I_t^Z- \bigr)\nonumber
\end{eqnarray}
for $t \ge0$. Inserting (\ref{28}) into (\ref{27}) and combining
it with (\ref{26}), we find that $V_1 = 1 - V_2$ for any $c>0$, and
this completes the proof.
\end{pf}

It follows from the previous proof that a stopping time $\tau$ of
$Z$ is optimal in~(\ref{23}) if and only if it is optimal in
(\ref{25}). To gain tractability when solving the optimal stopping
problem (\ref{25}) we will assume that the distribution function
$F$ of $\ell$ is strictly increasing and twice continuously
differentiable. Then $F(Z)$ defines a regular diffusion process with
values in $(0,1)$ and to gain symmetry and extend the state space to
$(-1,1)$, we will rescale $Z$ differently by setting
%
%e2.9 #&#
\begin{equation}
\label{29} X = 2 F(Z) - 1.
\end{equation}
Then $X$ is a regular diffusion process starting at $2F(0) - 1$ and
solving
%
%e2.10 #&#
\begin{equation}
\label{210} dX_t = \mu(X_t) \,dt +
\sigma(X_t) \,dB_t,
\end{equation}
where the drift $\mu$ and the diffusion coefficient $\sigma$ are
given by
%
%e2.11 #&#
%e2.12 #&#
\begin{eqnarray}
\label{211} \mu(x) &=& \bigl(2aF' + b^2
F'' \bigr) \biggl(F^{-1} \biggl(
\frac
{x+1}{2} \biggr) \biggr),
\\
\label{212} \sigma(x) &=& \bigl(2bF' \bigr)
\biggl(F^{-1} \biggl( \frac{x+1}{2} \biggr) \biggr)
\end{eqnarray}
for $x \in(-1,1)$ as is easily verified by It\^o's formula. Setting
%
%e2.13 #&#
\begin{equation}
\label{213} I_t = \inf_{0 \le s \le t} X_s
\quad\mbox{and}\quad S_t = \sup_{0
\le s
\le t}
X_s
\end{equation}
for $t \ge0$, we see that problem (\ref{25}) is equivalent to
the optimal stopping problem
%
%e2.14 #&#
\begin{equation}
\label{214} V = \sup_\tau\mathsf E \biggl[
S_\tau - I_\tau- c \int_0^\tau
( S_t - I_t ) \,dt \biggr],
\end{equation}
where the infimum is taken over all stopping times $\tau$ of $X$.
Note that $V = 2 V_2 = 2(1 - V_1)$, and there is a simple one-to-one
correspondence between the optimal stopping times in (\ref{214})
and (\ref{25}) due to (\ref{29}). We will therefore proceed by
studying problem (\ref{214}).

For future reference let us note that the infinitesimal generator of
$X$ equals
%
%e2.15 #&#
\begin{equation}
\label{215} \mathbb{L}_X = \mu(x) \frac{\partial}{\partial x} +
\frac{\sigma^2(x)}{2} \frac{\partial^2}{\partial x^2}
\end{equation}
and the scale function of $X$ is given by
%
%e2.16 #&#
\begin{equation}
\label{216} L(x) = \int_0^x \exp \biggl( -
\int_0^y \frac{\mu(z)}{(\sigma^2
/2)(z)} \,dz \biggr) \,dy
\end{equation}
for $x \in(-1,1)$. Throughout we denote $\rho_a = \inf\{ t \ge
0 \vert X_t=a \}$ and set $\rho_{a,b} = \rho_a \wedge\rho_b$
for $a<b$ in $(-1,1)$. Denoting by $\mathsf P_{ x}$ the probability
measure under which the process $X$ starts at $x$, it is well known
that
%
%e2.17 #&#
\begin{equation}
\label{217} \qquad\quad\mathsf P_{ x} ( X_{\rho_{a,b}}=a ) =
\frac{L(b) - L(x)}{L(b)
- L(a)} \quad\mbox{and}\quad\mathsf P_{ x} (
X_{\rho_{a,b}}=b ) = \frac{L(x) - L(a)} {L(b) - L(a)}
\end{equation}
for $a \le x \le b$ in $(-1,1)$. The speed measure of $X$ is given
by
%
%e2.18 #&#
\begin{equation}
\label{218} m(dx) = \frac{dx}{L'(x) (\sigma^2 /2)(x)}
\end{equation}
and the Green function of $X$ is given by
%
%e2.19 #&#
\begin{eqnarray}
\label{219} G_{a,b}(x,y) &=& \frac{(L(b) - L(y)) (L(x) - L(a))}{L(b)
- L(a)}\qquad\mbox{if } a \le
x \le y \le b
\nonumber\\[-8pt]\\[-8pt]
&=& \frac{(L(b) - L(x)) (L(y) - L(a))}{L(b) - L(a)}\qquad\mbox{if } a \le y \le x \le b.\nonumber
\end{eqnarray}
If $f\dvtx  (-1,1) \rightarrow\mathbb{R}$ is a measurable function, then
it is
well known that
%
%e2.20 #&#
\begin{equation}
\label{220} \mathsf E _x \int_0^{\rho_{a,b}}
f(X_t) \,dt = \int_a^b f(y)
G_{a,b}(x,y) m(dy)
\end{equation}
for $a \le x \le b$ in $(-1,1)$. This identity holds in the sense
that if one of the integrals exists, so does the other one, and they
are equal.

%s3 #&#
\section{Optimal stopping of the range process}\label{sec3}
%%%%%%%%%%%%%%%%%%%%%%%%%%%%%%%%%%%%%%%%%%%%%%%

It was shown in the previous section that the quickest detection
problem (\ref{23}) is equivalent to the optimal stopping problem
(\ref{214}). The purpose of this section is to present the solution
to the latter problem in somewhat greater generality. Using the fact
that the two problems are equivalent, this also leads to the solution
of the former problem.

Let $X=(X_t)_{t \ge0}$ be a one-dimensional diffusion process
solving
%
%e3.1 #&#
\begin{equation}
\label{31} dX_t = \mu(X_t) \,dt +
\sigma(X_t) \,dB_t,
\end{equation}
where the drift $\mu$ and the diffusion coefficient $\sigma>0$ are
continuous functions and $B=(B_t)_{t \ge0}$ is a standard Brownian
motion. To meet a sufficient condition used in the proof below, we
will also assume that $\sigma^2$ is (locally) Lipschitz. We will
further assume that the state space of $X$ equals $(-1,1)$ as in the
previous section; however, this hypothesis is not essential; see
Remark~\ref{rem4} below. By $\mathsf P_{ x}$ we denote the probability measure
under which $X$ starts at $x \in(-1,1)$. For $i \le x \le s$ in
$(-1,1)$ we set
%
%e3.2 #&#
\begin{equation}
\label{32} I_t = i \wedge\inf_{0 \le s \le t}
X_s \quad\mbox{and}\quad S_t = s \vee\sup
_{0 \le s \le t} X_s
\end{equation}
for $t \ge0$. These transformations enable the three-dimensional
Markov process $(I,X,S)$ to start at $(i,x,s)$ under $\mathsf P_{ x}$, and
we will denote the resulting probability measure on the canonical
space by $\mathsf P_{ i,x,s}$. Thus under $\mathsf P_{ i,x,s}$ the canonical
process $(I,X,S)$ starts at $(i,x,s)$. The range process $R$ of $X$
is defined by
%
%e3.3 #&#
\begin{equation}
\label{33} R_t = S_t - I_t
\end{equation}
for $t \ge0$. In this section we consider the optimal stopping
problem
%
%e3.4 #&#
\begin{equation}
\label{34} V(i,x,s) = \sup_\tau\mathsf E _{i,x,s}
\biggl[ R_\tau- \int_0^\tau
c(I_t,X_t,S_t) \,dt \biggr]
\end{equation}
for $i \le x \le s$ in $(-1,1)$ where the supremum is taken over all
stopping times $\tau$ of~$X$.

Regarding the cost function $c$ in (\ref{34}) we will assume that
(i) $i \mapsto c(i,x,s)$ is decreasing and $s \mapsto c(i,x,s)$ is
increasing with $c(i,x,s)>0$ for $i \le x \le s$ in~$(-1,1)$. These
conditions have a natural interpretation in the sense that any new
increase in gain (when $X$ reaches either $S$ or $I$) is followed by
a proportional increase in cost. To gain existence and tractability
we will also assume that (ii)~$(i,x,s) \mapsto c(i,x,s)$ is
continuous, $x \mapsto c(i,x,s)$ is (locally) Lipschitz, $(i,s)
\mapsto c(i,x,s)$ is continuously differentiable. To gain
monotonicity and joint continuity we will further assume that (iii)
$i \mapsto\frac{\partial c}{\partial s}(i,x,s)$ and $s \mapsto
\frac{\partial c}{\partial i}(i,x,s)$ are increasing and (locally)
Lipschitz. Note that conditions (i)--(iii) are satisfied for
$c(i,x,s) = c(s - i) > 0$ when $c$ is increasing concave and
continuously differentiable with $c'$ (locally) Lipschitz. Note also
that conditions (i)--(iii) are satisfied for $c(i,x,s) = c_2(s)
- c_1(i) > 0$ when $c_1$ and $c_2$ are increasing and continuously
differentiable functions. Note finally that conditions (i)--(iii)
are satisfied for $c(i,x,s)=c(x)>0$ when $c$ is (locally) Lipschitz
(in this case $f_*$ and $g_*$ below are no longer surfaces but
curves as functions of $i$ and $s$, respectively).

For any $s$ given and fixed we will refer to $d^s = \{ (i,x)
\vert i = x \le s \}$ as the \emph{lower diagonal} in the state
space, and for any $i$ given and fixed we will refer to $d_i = \{
(x,s) \vert x = s \ge i \}$ as the \emph{upper diagonal} in
the state space. We will say that a function $f$ \emph{stays
strictly above the lower diagonal} $d^s$ if $f(i,s)>i$ for all
$i<s$, and we will say that a function $g$ \emph{stays strictly
below the upper diagonal} $d_i$ if $g(i,s)<s$ for all $s>i$.

The main result of the paper may now be stated as follows.

%th1 #&#
\begin{teo}\label{teo1}
Under the hypotheses on $X$ and $c$ stated
above, the optimal stopping time in problem (\ref{34}) is given
by
%
%e3.5 #&#
\begin{equation}
\label{35} \tau_* = \inf \bigl\{ t \ge0 \vert f_*(I_t,S_t)
\le X_t \le g_*(I_t,S_t) \bigr\},
\end{equation}
where the surfaces $f_*$ and $g_*$ can be characterised as the
minimal and maximal solutions to
%
%e3.6 #&#
%e3.7 #&#
\begin{eqnarray}
\frac{\partial f}{\partial i}(i,s) &=& \frac{(\sigma^2 /2)(f(i,s))
L'(f(i,s))}{c(i,f(i,s),s) [L(f(i,s)) - L(i)]}
\nonumber\\[-8pt]\label{36} \\[-8pt]
&&{}\times \biggl[ 1 - \int_i^{f(i,s)} \frac{\partial c}{\partial i}(i,y,s)
\frac{L(y) - L(i)}{(\sigma^2 /2)(y) L'(y)} \,dy \biggr],\nonumber
\\
\frac{\partial g}{\partial s}(i,s) &=& \frac{(\sigma^2
/2)(g(i,s)) L'(g(i,s))}{c(i,g(i,s),s) [L(s) - L(g(i,s))]}
\nonumber\\[-8pt]\label{37} \\[-8pt]
&&{}\times  \biggl[ 1 + \int
_{g(i,s)}^s \frac{\partial c}{\partial s}(i, y,s)
\frac{L(s) - L(y)}{ (\sigma^2 /2)(y) L'(y)} \,dy \biggr]\nonumber
\end{eqnarray}
staying strictly above the lower diagonal $d^s$ and strictly below
the upper diagonal $d_i$ for $i<s$ in $(-1,1)$, respectively.

Explicit formulae for the value function $V$ on the
continuation sets (\ref{310}) and~(\ref{311}) below are given by
(\ref{325}) and (\ref{332}) below for any cost function $c$
satisfying \textup{(i)--(iii)} above. Explicit formulae for the value function
$V$ on the continuation set (\ref{39}) below are given by
(\ref{364}) and (\ref{365}) below when \mbox{$c(i,s) = c_2(s) - c_1(i)
> 0$} where $c_1$ and $c_2$ are increasing and continuously
differentiable functions. Outside these sets the value function $V$
equals $s - i$ for $i<s$ in $(-1,1)$. The optimal surfaces $f_*$
and $g_*$ satisfy the additional properties
\mbox{(\ref{334})--(\ref{339})}.
\end{teo}

\begin{pf}
The optimal stopping problem (\ref{34}) is
three-dimensional and the underlying Markov process equals
$(I,X,S)$. It is evident from the structure of the gain function in
(\ref{34}) that the excursions of $X$ away from the running maximum
$S$ and the running minimum $I$ play a key role in the analysis of
the problem. A~possible way to visualise the dynamics of these
excursions is illustrated in Figure~\ref{fig1} below. Each excursion of $X$
at an upper level $s$ is mirror imaged with the excursion of $X$ at
a lower level $i$ and vice versa. When the excursion returns to the
upper diagonal, the process $(X,S)$ receives an infinitesimal push
upwards along the upper diagonal, and when the excursion returns to
the lower diagonal, the process $(I,X)$ receives an infinitesimal
push downwards along the lower diagonal.

%%%%%%%%%%%%%%%%%%%%%%%%%%%%%%%%%%%%%%%%%%%%%%%%%%%%%%%%%%%%%%%%%%%%%%%%%%%%%%%%%
%%% Figure 1 %%%
%%%%%%%%%%%%%%%%%%%%%%%%%%%%%%%%%%%%%%%%%%%%%%%%%%%%%%%%%%%%%%%%%%%%%%%%%%%%%%%%%

%f1 #&#
\begin{figure}

\includegraphics{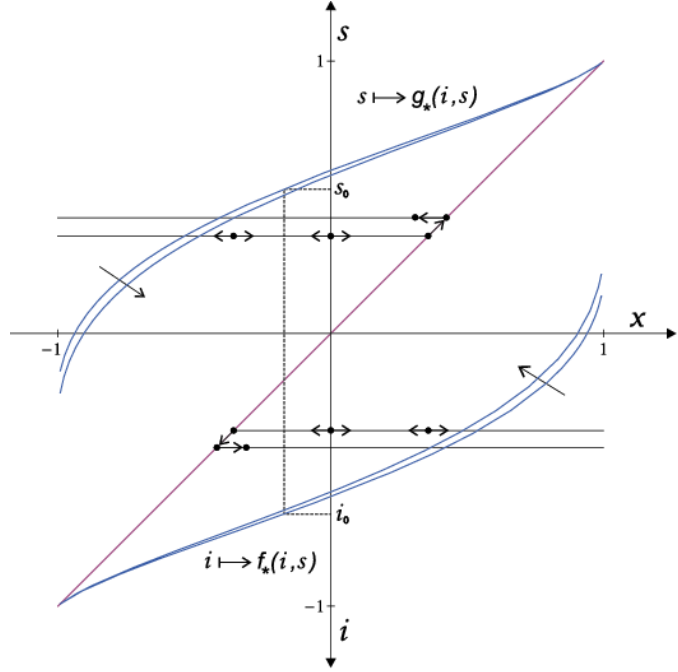}

\caption{Excursions of $X$ away from the running minimum
$I$ and the running maximum $S$ combined with the dynamics of the
optimal stopping surfaces $f_*$ and $g_*$:
\textup{(i)}~return of $X$ to the lower diagonal causes $I$ to go down and forces $g_*$ to go up;
\textup{(ii)}~return of $X$ to the upper diagonal causes $S$ to go up and forces $f_*$ to go down;
\textup{(iii)}~even if $X$ goes above $f_*$ it may not be optimal to stop unless $X$ is below $g_*$;
\textup{(iv)}~even if $X$ goes below $g_*$ it may not be optimal to stop unless $X$ is above $f_*$.
The (movable) dotted vertical line marks the borderline levels $i_0$
and $s_0$ below and above which it is optimal to stop.}\label{fig1}
\end{figure}

%%%%%%%%%%%%%%%%%%%%%%%%%%%%%%%%%%%%%%%%%%%%%%%%%%%%%%%%%%%%%%%%%%%%%%%%%%%%%%%%%
%%% End Figure 1 %%%
%%%%%%%%%%%%%%%%%%%%%%%%%%%%%%%%%%%%%%%%%%%%%%%%%%%%%%%%%%%%%%%%%%%%%%%%%%%%%%%%%

An important initial observation is that the process $(I,X,S)$ can
never be optimally stopped at the upper or lower diagonal. The
analogous phenomenon is known to hold for optimal stopping of the
maximum process (see \cite{Pe-1}, Proposition~2.1) and the same
arguments extend to the present case without major changes. Before
we formalise this in the first step below let us recall that general
theory of optimal stopping for Markov processes (see \cite{PS}, Chapter~1) implies that the continuation set in the problem (\ref{34})
equals $C = \{ (i,x,s) \vert V(i,x,s) > s - i \}$ and the
stopping set equals $D = \{ (i,x,s) \vert V(i,x,s) = s - i
\}$. It means that the first entry time of $(I,X,S)$ into $D$ is
optimal in problem (\ref{34}). To determine the sets $C$ and
$D$ we will begin by formalising the initial observation above.
\begin{longlist}[(3)]
\item[(1)] \emph{The upper and lower diagonal $d_i$ and $d^s$ are always
contained in $C$}. For this, take any $(s,s) \in d_i$ and consider
$\rho_{l_n,r_n} = \inf\{ t \ge0 \vert X_t \notin
(l_n,r_n) \}$ under $\mathsf P_{ i,s,s}$ with $l_n = s - 1/n$ and $r_n
= s + 1/n$ for $n \ge1$. Then (\ref{218})--(\ref{220}) imply that
$\mathsf E _{i,s,s} R_{\rho_{l_n,r_n}} \ge s - i + K/n$ and $\mathsf E _{i,s,s}
\int_0^{\rho_{l_n,r_n}} c(I_t,X_t,S_t) \,dt \le K/n^2$ for all $n
\ge1$ with some positive constant $K$ (see the proof of Proposition~2.1 in \cite{Pe-1} for details). Taking $n \ge1$ large enough (to
exploit the difference in the rates of the bounds) we see that
$(i,s,s)$ belongs to $C$. In exactly the same way one sees that if
$(i,i) \in d^s$ then $(i,i,s)$ belongs to $C$. This establishes the
initial claim.
\end{longlist}
\begin{longlist}[(3)]
\item[(2)] \emph{Optimal stopping surfaces}. Assume now that the process
$(I,X,S)$ starts at $(i,x,s)$, and consider the excursion of $X$ away
from the running maximum $s$ with $i$ given and fixed. In view of
the fact that it is never optimal to stop at the upper diagonal
$d_i$, and due to the existence of a strictly positive cost which is
proportional to the duration of time in (\ref{34}), we see that it
is plausible to expect that there exists a point $g(i,s)$ (depending
on both $i$ and $s$) at/below which the process $X$ should be
stopped (should $i$ remain constant). In exactly the same way, if we
consider the excursion of $X$ away from the running minimum $i$ with
$s$ given and fixed, we see that it is plausible to expect that
there exists a point $f(i,s)$ (depending on both $i$~and~$s$)
at/above which the process $X$ should be stopped (should $s$ remain
constant).

The first complication in this reasoning comes from the fact that
neither $i$ nor $s$ need to remain constant during the excursion of
$X$ away from the running maximum $s$ or the running minimum $i$,
respectively. We will handle this difficulty implicitly by noting
that if $I$ is to decrease from $i$ downwards, then this will
increase the rate of the cost in (\ref{34}) which in turn will move
the boundary point $g(i,s)$ upwards [it means that $i \mapsto
g(i,s)$ is decreasing], and similarly if $S$ is to increase from $s$
upwards then this will increase the rate of the cost in (\ref{34})
which in turn will move the boundary point $f(i,s)$ downwards [it
means that $s \mapsto f(i,s)$ is decreasing]. To visualise these
movements see Figure~\ref{fig1} above. Changes in either $I$ or $S$ therefore
contribute to resetting $i$ and $s$ to new levels and starting from
there afresh with the boundary points $f(i,s)$ and $g(i,s)$
adjusted. For these reasons it is not entirely surprising that the
first complication will resolve itself after we describe the
structure of the optimal surfaces $f$ and $g$ in fuller detail
below.

The second complication comes from the fact that even if $X$ is
at/below $g(i,s)$ and normally (when $i$ would not change) it would
be optimal to stop, it may be that $X$ is still below $f(i,s)$ and
therefore the proximity of the lower diagonal $d^s$ may be a valid
incentive to continue. This incentive itself is further complicated
by the fact that it may lead to a decrease of $i$ and therefore the
rate of the cost in~(\ref{34}) will also increase (as addressed in
the first complication above). Likewise, even if $X$ is at/above
$f(i,s)$ and normally (when $s$ would not change) it would be
optimal to stop, it may be that $X$ is still above $g(i,s)$ and
therefore the proximity of the upper diagonal $d_i$ may be a valid
incentive to continue. This incentive itself is further complicated
by the fact that it may lead to an increase of $s$ and therefore the
rate of the cost in (\ref{34}) will also increase (as addressed in
the first complication above).

Neither of these complications appear in the optimal stopping of the
maximum process where $g$ depends only on $s$ (see \cite{Pe-1} and
the references therein), and our strategy in tackling the problem
will be to extend the maximality principle \cite{Pe-1} from the
two-dimensional setting of the process $(X,S)$ and the optimal
stopping curves to the three-dimensional setting of the process
$(I,X,S)$ and the optimal stopping surfaces. This will enable us to
resolve the second complication using the existence of the so-called
``bad--good'' solutions (those hitting the upper or lower diagonal)
which in turn will provide novel insights into the
maximality/minimality principle in the three dimensions as will be
seen below.
\end{longlist}
\begin{longlist}[(3)]
\item[(3)] \emph{Free-boundary problem}. Previous considerations suggest to
seek the solution to (\ref{34}) as the following stopping time:
%
%e3.8 #&#
\begin{equation}
\label{38} \tau_{f,g} = \inf \bigl\{ t \ge0 \vert
f(I_t,S_t) \le X_t \le g(I_t,S_t)
\bigr\},
\end{equation}
where the surfaces $f$ and $g$ are to be found. The continuation set
$C_{f,g}$ splits into
%
%e3.9 #&#
%e3.10 #&#
%e3.11 #&#
\begin{eqnarray}
\label{39} C_{f,g}^0 &=& \bigl\{ (i,x,s) \vert f(i,s) >
g(i,s) \bigr\},
\\
\label{310} C_{f,g}^- &=& \bigl\{ (i,x,s) \vert i \le x < f(i,s) \le
g(i,s) \bigr\},
\\
\label{311} C_{f,g}^+ &=& \bigl\{ (i,x,s) \vert f(i,s) \le g(i,s) < x
\le s \bigr\}
\end{eqnarray}
and we have $C_{f,g} = C_{f,g}^0 \cup C_{f,g}^- \cup C_{f,g}^+$. To
compute the value function $V$ and determine the optimal surfaces
$f$ and $g$, we are led to formulate the free-boundary problem
%
%e3.12 #&#
%e3.13 #&#
%e3.14 #&#
%e3.15 #&#
%e3.16 #&#
%e3.17 #&#
%e3.18 #&#
\begin{eqnarray}
\label{312} (\mathbb{L}_X V) (i,x,s) &=& c(i,x,s)\qquad\mbox{for
} (i,x,s) \in C_{f,g},
\\
\label{313} V_i'(i,x,s) \vert_{x=i+} &=& 0
\qquad\mbox {(normal reflection)},
\\
\label{314} V_s'(i,x,s) \vert_{x=s-} &=& 0
\qquad\mbox {(normal reflection)},
\\
\label{315} V(i,x,s) \vert_{x=f(i,s)-} &=& s - i\qquad\mbox {for } f(i,s)
\le g(i,s),
\\
\label{316} V(i,x,s) \vert_{x=g(i,s)+} &=& s - i\qquad\mbox {for } f(i,s)
\le g(i,s),
\\
\label{317} V_x'(i,x,s) \vert_{x=f(i,s)-} &=& 0
\qquad\mbox {for } f(i,s) \le g(i,s)\qquad\mbox{(smooth fit)},
\\
\label{318} \qquad V_x'(i,x,s) \vert_{x=g(i,s)+} &=& 0
\qquad\mbox {for } f(i,s) \le g(i,s)\qquad\mbox{(smooth fit)},
\end{eqnarray}
where $\mathbb{L}_X$ is the infinitesimal generator of $X$ given in
(\ref{215}) above. For the rationale and further details regarding
free-boundary problems of this kind, we refer to~\cite{PS}, Section~13, and the references therein; we note in addition that the
conditions of normal reflection (\ref{313}) and (\ref{314}) date
back to \cite{GSG}.
\end{longlist}
\begin{longlist}[(3)]
\item[(4)] \emph{Nonlinear differential equations}. To tackle the
free-boundary problem (\ref{312})--(\ref{318}), consider the
resulting function
%
%e3.19 #&#
\begin{equation}
\label{319} V_{f,g}(i,x,s) = \mathsf E _{i,x,s} \biggl[
R_{\tau_{f,g}} - \int_0^{
\tau_{f,g}}
c(I_t,X_t,S_t) \,dt \biggr]
\end{equation}
for $i \le x \le s$ in $(-1,1)$ upon assuming that $\mathsf E _{i,x,s}
\tau_{f,g} < \infty$ with candidate surfaces $f$ and $g$ to be
specified below. Suppose that $f(i,s) \le s$ and consider
$\rho_{i,f(i,s)} = \inf\{ t \ge0 \vert X_t \notin
(i,f(i,s)) \}$ under $\mathsf P_{i,x,s}$ with $i<x<f(i,s)$ given and
fixed. Applying the strong Markov property of $(I,X,S)$ at
$\rho_{i,f(i,s)}$ and using (\ref{217})--(\ref{220}) we find that
%
%e3.20 #&#
\begin{eqnarray}
\label{320} V_{f,g}(i,x,s) &=& (s - i) \frac{L(x) - L(i)}{L(f(i,s)) - L(i)}\nonumber
\\
&&{} + V_{f,g}(i,i,s) \frac{L(f(i,s)) - L(x)}{L(f(i,s)) - L(i)}
\\
&&{} - \int_i^{f(i,s)} c(i,y,s)
G_{i,f(i,s)}(x,y) m(dy).\nonumber
\end{eqnarray}
It follows from (\ref{320}) that
%
%e3.21 #&#
\begin{eqnarray}
\label{321} V_{f,g}(i,i,s) &=& s - i\nonumber
\\
&&{}  + \frac{L(f(i,s)) - L(i)}{L(f(i,s)) - L(x)}
\biggl[V_{f,g}(i,x,s) - (s - i)
\\
&&\hspace*{99pt}{}+ \int_i^{f(i,s)} c(i,y,s)
G_{i,f(i,s)}(x,y) m(dy) \biggr].\hspace*{-10pt}\nonumber
\end{eqnarray}
Dividing and multiplying through by $x - f(i,s)$ we find using
(\ref{317}) that
%
%e3.22 #&#
\begin{eqnarray}
\label{322} && \lim_{x \uparrow f(i,s)} \frac{V_{f,g}(i,x,s) - (s - i)}{L(f(i,s))
- L(x)}
\nonumber\\[-8pt]\\[-8pt]
&&\qquad =  -\frac{1}{L'(f(i,s))} \frac{\partial V_{f,g}}{\partial
x} (i,x,s) \bigg\vert_{x = f(i,s)-}
 = 0\nonumber
\end{eqnarray}
for $f(i,s) \le g(i,s)$. It is easily seen by (\ref{219}) that
%
%e3.23 #&#
\begin{eqnarray}
\label{323} &&\lim_{x \uparrow f(i,s)} \frac{L(f(i,s)) - L(i)}{L(f(i,s)) - L(x)} \int
_i^{f(i,s)} c(i,y,s) G_{i,f(i,s)}(x,y) m(dy)
\nonumber\\[-8pt]\\[-8pt]
&&\qquad = \int_i^{f(i,s)} c(i,y,s) \bigl[L(y) -
L(i) \bigr] m(dy).\nonumber
\end{eqnarray}
Combining (\ref{321})--(\ref{323}) we find that
%
%e3.24 #&#
\begin{equation}
\label{324} V_{f,g}(i,i,s) = s - i + \int_i^{f(i,s)}
c(i,y,s) \bigl[L(y) - L(i) \bigr] m(dy)
\end{equation}
for $f(i,s) \le g(i,s)$. Inserting this back into (\ref{320}) and
using (\ref{219}) and (\ref{220}) we conclude that
%
%e3.25 #&#
\begin{equation}
\label{325} V_{f,g}(i,x,s) = s - i + \int_x^{f(i,s)}
c(i,y,s) \bigl[L(y) - L(x) \bigr] m(dy)
\end{equation}
for $x \le f(i,s) \le g(i,s)$. Finally, using (\ref{313}) we find
that
%
%e3.26 #&#
\begin{eqnarray}
\label{326} \frac{\partial f}{\partial i}(i,s) &=& \frac{(\sigma^2 /2)(f(i,s))
L'(f(i,s))}{c(i,f(i,s),s) [L(f(i,s)) - L(i)]}
\nonumber\\[-8pt]\\[-8pt]
&&{}\times  \biggl[ 1 - \int_i^{f(i,s)} \frac{\partial c}{\partial i}(i,y,s) \bigl[L(y) -
L(i) \bigr] m(dy) \biggr]\nonumber
\end{eqnarray}
for $f(i,s) \le g(i,s)$. By (\ref{218}) we see that (\ref{326})
coincides with (\ref{36}) above.

Similarly, suppose that $g(i,s) \ge i$ and consider $\rho_{g(i,s),s}
= \inf\{ t \ge0 \vert X_t \notin(g(i,s),s) \}$ under
$\mathsf P_{i,x,s}$ with $g(i,s)<x<s$ given and fixed. Applying the strong
Markov property of $(I,X,S)$ at $\rho_{g(i,s),s}$ and using
(\ref{217})--(\ref{220}) we find that
%
%e3.27 #&#
\begin{eqnarray}
\label{327} V_{f,g}(i,x,s) &=& (s - i) \frac{L(s) - L(x)}{L(s) - L(g(i,s))}\nonumber
\\
&&{}  +V_{f,g}(i,s,s) \frac{L(x) - L(g(i,s))}{L(s) - L(g(i,s))}
\\
&&{} - \int_{g(i,s)}^s c(i,y,s)
G_{g(i,s),s}(x,y) m(dy).\nonumber
\end{eqnarray}
It follows from (\ref{327}) that
%
%e3.28 #&#
\begin{eqnarray}
\label{328} V_{f,g}(i,s,s) &=& s - i\nonumber
\\
&&{}  + \frac{L(s) - L(g(i,s))}{L(x) - L(g(i,s))} \biggl[
V_{f,g}(i,x,s) - (s - i)
\\
&&\hspace*{98pt}{}+ \int_{g(i,s)}^s c(i,y,s)
G_{g(i,s),s}(x,y) m(dy) \biggr].\hspace*{-10pt}\nonumber
\end{eqnarray}
Dividing and multiplying through by $x - g(i,s)$ we find using
(\ref{318}) that
%
%e3.29 #&#
\begin{eqnarray}\label{329}
&& \lim_{x \downarrow g(i,s)} \frac{V_{f,g}(i,x,s) - (s - i)}{L(x)
- L(g(i,s))}
\nonumber\\[-8pt]\\[-8pt]
&&\qquad =\frac{1}{L'(g(i,s))} \frac{\partial V_{f,g}}{
\partial x} (i,x,s) \bigg\vert_{x = g(i,s)+} = 0\nonumber
\end{eqnarray}
for $g(i,s) \ge f(i,s)$. It is easily seen by (\ref{219}) that
%
%e3.30 #&#
\begin{eqnarray}
\label{330} &&\lim_{x \downarrow g(i,s)} \frac{L(s) - L(g(i,s))}{L(x) - L(g(i,s))} \int
_{g(i,s)}^s c(i,y,s) G_{g(i,s),s}(x,y) m(dy)
\nonumber\\[-8pt]\\[-8pt]
&&\qquad = \int_{g(i,s)}^s c(i,y,s) \bigl[L(s) -
L(y) \bigr] m(dy).\nonumber
\end{eqnarray}
Combining (\ref{328})--(\ref{330}) we find that
%
%e3.31 #&#
\begin{equation}
\label{331} V_{f,g}(i,s,s) = s - i + \int_{g(i,s)}^s
c(i,y,s) \bigl[L(s) - L(y) \bigr] m(dy)
\end{equation}
for $g(i,s) \ge f(i,s)$. Inserting this back into (\ref{327}) and
using (\ref{219}) and (\ref{220}) we conclude that
%
%e3.32 #&#
\begin{equation}
\label{332} V_{f,g}(i,x,s) = s - i + \int_{g(i,s)}^x
c(i,y,s) \bigl[L(x) - L(y) \bigr] m(dy)
\end{equation}
for $x \ge g(i,s) \ge f(i,s)$. Finally, using (\ref{314}) we find
that
%
%e3.33 #&#
\begin{eqnarray}\label{333}
\frac{\partial g}{\partial s}(i,s) &=& \frac{(\sigma^2 /2)(g(i,s))
L'(g(i,s))}{c(i,g(i,s),s) [L(s) - L(g(i,s))]}
\nonumber\\[-8pt]\\[-8pt]
&&{}\times  \biggl[ 1 + \int
_{g(i,s)}^s \frac{\partial c}{\partial s}(i,y,s) \bigl[L(s) - L(y)
\bigr] m(dy) \biggr]\nonumber
\end{eqnarray}
for $g(i,s) \ge f(i,s)$. By (\ref{218}) we see that (\ref{333})
coincides with (\ref{37}) above.

Summarising the preceding considerations we can conclude that to
each pair of the candidate surfaces\vspace*{1pt} $f$ and $g$ solving (\ref{36})
and (\ref{37}) there corresponds the function
(\ref{325}) and (\ref{332}) on $C_{f,g}^- \cup C_{f,g}^+$ solving the
free-boundary problem \mbox{(\ref{312})--(\ref{318})} on $C_{f,g}^- \cup
C_{f,g}^+$ (this\vspace*{-1pt} can be verified by direct differentiation) and
admitting the probabilistic representation\vspace*{1.5pt} (\ref{319}) on
$C_{f,g}^- \cup C_{f,g}^+$ associated with the stopping time
(\ref{38}) when the latter has finite expectation [this will be
formally proved for the surfaces of interest in (\ref{374}) and
(\ref{375}) below].

The central question becomes how to select the optimal surfaces $f$
and $g$ among all admissible candidates solving (\ref{36}) and
(\ref{37}). We will answer this question by invoking the
superharmonic characterisation of the value function (see
\cite{PS}, Chapter~1) for the four-dimensional Markov process
$(I,X,S,A)$ where $A_t = \int_0^t c(I_s,X_s,S_s) \,ds$ for $t \ge
0$. Fuller details of this argument will become clearer as we
progress below.
\end{longlist}
\begin{longlist}[(3)]
\item[(5)] \emph{The minimal and maximal solution}. Motivated by the
previous question we note from (\ref{325}) and (\ref{332}) that $f
\mapsto V_{f,g}$ is increasing and $g \mapsto V_{f,g}$ is
decreasing. Recalling also that it is not optimal to stop at the
upper or lower diagonal, this motivates us to select solutions to
(\ref{36}) and (\ref{37}) as far as possible from the upper and
lower diagonal, respectively [respecting also the meaning of
(\ref{38}) in (\ref{319}) as well as the meaning of (\ref{319})
itself]. In the former case this means as small as possible below
the upper diagonal, and in the latter case it means as large as
possible above the lower diagonal. We ought to recall, however, that
stopping time~(\ref{38}) needs to have finite expectation, and
this will put a natural constraint on how small and large these
solutions can be (this is a subtle point in the background of the
argument).

%%%%%%%%%%%%%%%%%%%%%%%%%%%%%%%%%%%%%%%%%%%%%%%%%%%%%%%%%%%%%%%%%%%%%%%%%%%%%%%%%
%%% Figure 2 %%%
%%%%%%%%%%%%%%%%%%%%%%%%%%%%%%%%%%%%%%%%%%%%%%%%%%%%%%%%%%%%%%%%%%%%%%%%%%%%%%%%%

%f2 #&#
\begin{figure}

\includegraphics{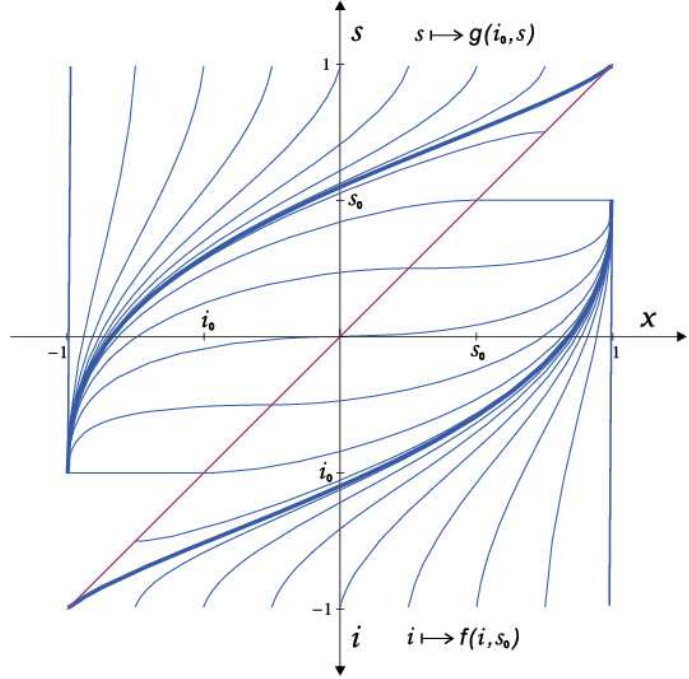}

\caption{Smooth-fit solutions $i \mapsto f(i,s_0)$ and $s
\mapsto g(i_0,s)$ to differential equations (\protect\ref{36}) and
(\protect\ref{37}) for fixed $s_0$ and $i_0$, respectively. The minimal
solution staying strictly above the lower diagonal (bold $f$ line)
and the maximal solution staying strictly below the upper diagonal
(bold $g$ line) are sections of the optimal stopping surfaces, respectively.}\label{fig2}
\end{figure}

%%%%%%%%%%%%%%%%%%%%%%%%%%%%%%%%%%%%%%%%%%%%%%%%%%%%%%%%%%%%%%%%%%%%%%%%%%%%%%%%%
%%% End Figure 2 %%%
%%%%%%%%%%%%%%%%%%%%%%%%%%%%%%%%%%%%%%%%%%%%%%%%%%%%%%%%%%%%%%%%%%%%%%%%%%%%%%%%%

To address the existence and uniqueness of solutions to these
equations, denote the right-hand side of~(\ref{36}) by
$\Phi(i,s,f(i,s))$ and denote the right-hand side of~(\ref{37}) by
$\Psi(i,s,g(i,s))$. From general theory of nonlinear differential
equations we know that if the direction fields $(i,f) \mapsto
\Phi(i,s,f)$ and $(s,g) \mapsto\Psi(i,s,g)$ are (locally)
continuous and (locally) Lipschitz in the second variable, then
equations (\ref{36}) and (\ref{37}) admit (locally) unique
solutions. In particular, recalling that $(i,x,s) \mapsto c(i,x,s)$
is continuous we see from the structure of $\Phi$ and $\Psi$ that
equations (\ref{36}) and (\ref{37}) admit (locally) unique
solutions since $x \mapsto\sigma^2(x)$ and $x \mapsto c(i,x,s)$ are
(locally) Lipschitz.

%%%%%%%%%%%%%%%%%%%%%%%%%%%%%%%%%%%%%%%%%%%%%%%%%%%%%%%%%%%%%%%%%%%%%%%%%%%%%%%%%
%%% Figure 3 %%%
%%%%%%%%%%%%%%%%%%%%%%%%%%%%%%%%%%%%%%%%%%%%%%%%%%%%%%%%%%%%%%%%%%%%%%%%%%%%%%%%%

%f3 #&#
\begin{figure}[b]

\includegraphics{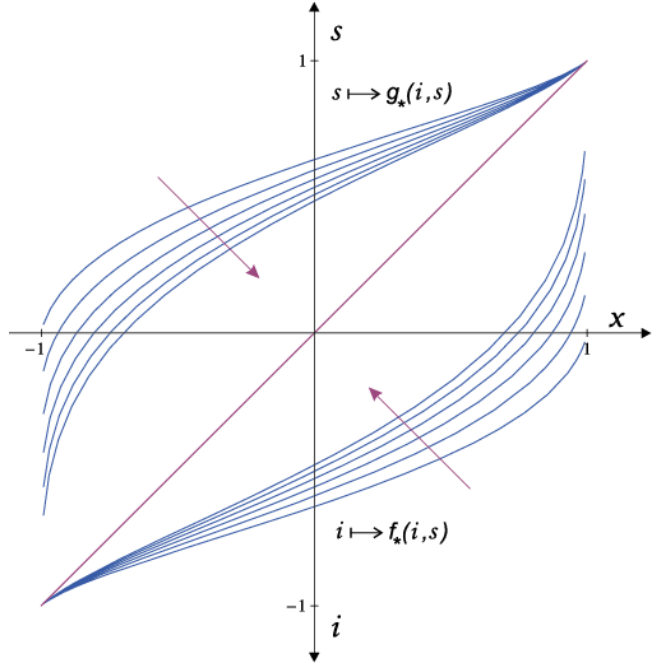}

\caption{Movement and shape of sections $i \mapsto
f_*(i,s)$ and $s \mapsto g_*(i,s)$ of the optimal surfaces $f_*$ and
$g_*$ as the running maximum $s$ increases and the running minimum
$i$ decreases, respectively.}\label{fig3}
\end{figure}

%%%%%%%%%%%%%%%%%%%%%%%%%%%%%%%%%%%%%%%%%%%%%%%%%%%%%%%%%%%%%%%%%%%%%%%%%%%%%%%%%
%%% End Figure 3 %%%
%%%%%%%%%%%%%%%%%%%%%%%%%%%%%%%%%%%%%%%%%%%%%%%%%%%%%%%%%%%%%%%%%%%%%%%%%%%%%%%%%

To construct the minimal solution to (\ref{36}) staying strictly
above the lower diagonal $d^s$, we can proceed as follows; see Figure~\ref{fig2} above. For any $i_n \in(-1,1)$ such that $i_n \downarrow-1$ as
$n \rightarrow\infty$ let $i \mapsto f_n(i,s)$ denote the solution
to (\ref{36}) such that $f_n(i_n,s)=i_n$ for $n \ge1$. Note that
each solution $i \mapsto f(i,s)$ to (\ref{36}) is singular at the
lower diagonal $d^s$ in the sense that $f_i'(i+,s) = +\infty$ for
$f(i+,s)=i$; however, passing to the equivalent equation for the
inverse of $i \mapsto f(i,s)$ [upon noting that each solution $i
\mapsto f(i,s)$ to (\ref{36}) is strictly increasing] we see that
this singularity gets removed; note that the inverse of $i
\mapsto f(i,s)$ has the derivative equal to zero at the lower
diagonal $d^s$. By the uniqueness of the solution we know that the
two curves $i \mapsto f_n(i,s)$ and $i \mapsto f_m(i,s)$ cannot
intersect for $n \ne m$, and hence we see that $(f_n)_{n \ge1}$ is
increasing. It follows therefore that $f_*:= \lim_{ n \rightarrow
\infty} f_n$ exists. Passing to an integral equation equivalent to
(\ref{36}) (or its inverse), it is easily verified that $i \mapsto
f_*(i,s)$ solves (\ref{36}) whenever strictly larger than $-1$.
This $f_*$ represents the minimal solution to (\ref{36}) staying
strictly above the lower diagonal. Since $i \mapsto c(i,x,s)$ is
decreasing we see from (\ref{36}) that
%
%e3.34 #&#
%e3.35 #&#
\begin{eqnarray}\label{334}
&& i \mapsto f_n(i,s)\mbox{ and } i \mapsto f_*(i,s)
\mbox{ are strictly increasing }
\nonumber\\[-10pt]\\[-10pt]
&&\mbox{with } f_*(-1+,s) = -1\nonumber
\end{eqnarray}
for $i<s$ in $(-1,1)$ and $n \ge1$. Note further that the increase
of $s \mapsto\frac{\partial c}{\partial i}(i,x,s)$ combined with
the increase of $s \mapsto c(i,x,s)$ implies that $s \mapsto
\Phi(i,s,f)$ is decreasing. Recalling that (\ref{36}) is being
solved forwards, this shows that
%
%e3.36 #&#
\begin{equation}
\label{335} s \mapsto f_n(i,s)\mbox{ and } s \mapsto f_*(i,s)\mbox{ are decreasing}
\end{equation}
for $i<s$ in $(-1,1)$ and $n \ge1$; see Figure~\ref{fig3} below. Moreover,
since $s \mapsto\frac{\partial c}{\partial i}(i,x,s)$ is (locally)
Lipschitz we see that $s \mapsto\Phi(i,s,f)$ is (locally) Lipschitz
from where we can easily deduce using Gronwall's inequality that
%
%e3.37 #&#
\begin{equation}
\label{336} (i,s) \mapsto f_n(i,s)\mbox{ and } (i,s) \mapsto f_*(i,s)\mbox{ are continuous}
\end{equation}
for $i<s$ in $(-1,1)$ and $n \ge1$. To simplify the notation we
will use the same symbol $f$ below to denote the minimal solution
$f_*$ unless stated otherwise.

To construct the maximal solution to (\ref{37}) staying strictly
below the upper diagonal $d_i$, we can proceed similarly; see Figure~\ref{fig2} above. For any $s_n \in(-1,1)$ such that $s_n \uparrow1$ as $n
\rightarrow\infty$ let $s \mapsto g_n(i,s)$ denote the solution to
(\ref{37}) such that $g_n(i,s_n)=s_n$ for $n \ge1$. Note that each
solution $s \mapsto g(i,s)$ to (\ref{37}) is singular at the upper
diagonal $d_i$ in the sense that $g_s'(i,s-) = +\infty$ for
$g(i,s-)=s$; however, passing to the equivalent equation for the
inverse of $s \mapsto g(i,s)$ [upon noting that each solution $s
\mapsto g(i,s)$ to (\ref{37}) is strictly increasing], we see that
this singularity gets removed; note that the inverse of $s
\mapsto g(i,s)$ has the derivative equal to zero at the upper
diagonal $d_i$. By the uniqueness of the solution we know that the
two curves $s \mapsto g_n(i,s)$ and $s \mapsto g_m(i,s)$ cannot
intersect for $n \ne m$, and hence we see that $(g_n)_{n \ge1}$ is
decreasing. It follows therefore that $g_*:= \lim_{ n \rightarrow
\infty} g_n$ exists. Passing to an integral equation equivalent to
(\ref{37}) (or its inverse) it is easily verified that $s \mapsto
g_*(i,s)$ solves (\ref{37}) whenever strictly smaller than $1$.
This $g_*$ represents the maximal solution to (\ref{37}) staying
strictly below the upper diagonal. Since $s \mapsto c(i,x,s)$ is
increasing we see from (\ref{37}) that
%
%e3.38 #&#
%e3.39 #&#
\begin{eqnarray}\label{337}
s \mapsto g_n(i,s)\mbox{ and } s \mapsto g_*(i,s)
\mbox{ are strictly increasing with } g_*(i,1-) = 1\hspace*{-35pt}
\end{eqnarray}
for $i<s$ in $(-1,1)$ and $n \ge1$. Note further that the increase
of $i \mapsto\frac{\partial c}{\partial s}(i,x,s)$ combined with
the decrease of $i \mapsto c(i,x,s)$ implies that $i \mapsto
\Psi(i,s,f)$ is increasing. Recalling that (\ref{37}) is being
solved backwards, this shows that
%
%e3.40 #&#
\begin{equation}
\label{338} i \mapsto g_n(i,s)\mbox{ and } i \mapsto g_*(i,s)\mbox{ are decreasing}
\end{equation}
for $i<s$ in $(-1,1)$ and $n \ge1$; see Figure~\ref{fig3} above. Moreover,
since $i \mapsto\frac{\partial c}{\partial s}(i,s)$ is (locally)
Lipschitz we see that $i \mapsto\Psi(i,s,f)$ is (locally) Lipschitz
from where we can easily deduce using Gronwall's inequality that
%
%e3.41 #&#
\begin{equation}
\label{339} (i,s) \mapsto g_n(i,s) \mbox{ and } (i,s) \mapsto g_*(i,s)
\mbox{ are continuous}
\end{equation}
for $i<s$ in $(-1,1)$ and $n \ge1$. To simplify the notation we
will use the same symbol $g$ below to denote the maximal solution
$g_*$ unless stated otherwise.

With the minimal and maximal solution $f$ and $g$ we can associate
the stopping time (\ref{38}) and the resulting function
(\ref{319}). Doing the same thing with $f_n$ and $g_n$ [noting that
the stopping time (\ref{38}) has finite expectation], the arguments
above show that (\ref{325}) and (\ref{332}) hold for $f_n$ and
$g_n$ for $n \ge1$. Passing in these expressions to the limit as $n
\rightarrow\infty$, we see that (\ref{325}) and (\ref{332}) remain
valid for the minimal and maximal solution $f$ and $g$. The claims
of the past two sentences will be formally verified in (\ref{374})
and (\ref{375}) below. This establishes closed-form expressions for
$V_{f,g}$ in terms of $f$ and $g$ on $C_{f,g}^+$ and $C_{f,g}^-$.
\end{longlist}
\begin{longlist}[(3)]
\item[(6)] \emph{Computing $V_{f,g}$ on $C_{f,g}^0$}. This calculation is
technically more complicated, and we will derive\vspace*{1pt} closed-form
expressions for $V_{f,g}$ in terms of $f$ and $g$ on $C_{f,g}^0$
when $c(i,s) = c_2(s) - c_1(i) > 0$ where $c_1$ and $c_2$ are
increasing and continuously differentiable functions. Note that the
latter decomposition is fulfilled in the setting in Section~\ref{sec2} above.
Note also that these closed-form expressions are not needed to
derive the optimality of $f$ and $g$ as it will be shown in the rest
of the proof below.

We begin\vspace*{1pt} by noting that $V_{f,g}$ needs to satisfy
(\ref{312})--(\ref{314}) on $C_{f,g}^0$; see Remark~\ref{rem2} below.
Recalling that a particular solution to $\mathbb{L}_X H=1$ is given by
%
%e3.42 #&#
\begin{equation}
\label{340} H(x) = \int_0^x \bigl[L(x) -
L(y) \bigr] m(dy),
\end{equation}
it follows from (\ref{312}) that
%
%e3.43 #&#
\begin{equation}
\label{341} V(i,x,s) = A(i,s) L(x) + B(i,s) + \bigl(c_2(s) -
c_1(i) \bigr) H(x)
\end{equation}
for some unknown functions $A$ and $B$ to be found. By (\ref{313})
and (\ref{314}) we find that
%
%e3.44 #&#
%e3.45 #&#
\begin{eqnarray}
\label{342} A_i'(i,s) L(i) + B_i'(i,s)
- c_1'(i) H(i)&=& 0,
\\
\label{343} A_s'(i,s) L(s) + B_s'(i,s)
+ c_2'(s) H(s) &=& 0.
\end{eqnarray}
Differentiating (\ref{342}) with respect to $s$ and (\ref{343})
with respect to $i$ (upon assuming that $A$ and $B$ are twice
continuously differentiable) it follows by subtracting the resulting
identities that $A_{is}''(i,s)=0$ and hence $B_{is}''(i,s)=0$ too.
This implies that
%
%e3.46 #&#
\begin{equation}
\label{344} A(i,s) = a_1(i) + a_2(s) \quad\mbox{and}\quad B(i,s) = b_1(i) + b_2(s)
\end{equation}
for some $a_i$ and $b_i$ to be found when $i=1,2$. Inserting this
back into (\ref{341})--(\ref{343}) we obtain
%
%e3.47 #&#
%e3.48 #&#
%e3.49 #&#
\begin{eqnarray}
\qquad && V(i,x,s)= \bigl(a_1(i) + a_2(s) \bigr) L(x)
\nonumber\\[-8pt]\label{345}\\[-8pt]
&&\hspace*{53pt} {}+ b_1(i)
+ b_2(s) + \bigl(c_2(s) -
c_1(i) \bigr) H(x),\nonumber
\\
\label{346}  && a_1'(i) L(i) + b_1'(i)
- c_1'(i) H(i) = 0,
\\
\label{347}  && a_2'(s) L(s) + b_2'(s)
+ c_2'(s) H(s) = 0
\end{eqnarray}
for $f(i,s) > g(i,s)$.

To determine $a_i$ and $b_i$ for $i=1,2$ recall that $V_{f,g}$ is
known at $C_{f,g}^-$ and $C_{f,g}^+$ so that it is also known at the
boundary between $C_{f,g}^0$ and $C_{f,g}^-$ and the boundary
between $C_{f,g}^0$ and $C_{f,g}^+$. This serves as a basic
motivation for the\vspace*{1pt} introduction of the following functions. Given
$(i,s)$ such that $f(i,s) > g(i,s)$ there exist unique $i(s)<i$ and
$s(i)>s$ such that
%
%e3.50 #&#
\begin{equation}
\label{348} f \bigl(i(s),s \bigr) = g \bigl(i(s),s \bigr) \quad\mbox{and}
\quad f \bigl(i,s(i) \bigr) = g \bigl(i,s(i) \bigr).
\end{equation}
The existence of $i(s)$ and $s(i)$ follows from the facts that $i
\mapsto f(i,s)$ and $s \mapsto g(i,s)$ are strictly increasing and
$s \mapsto f(i,s)$ and $i \mapsto g(i,s)$ are strictly decreasing; see
Figure~\ref{fig3} above. More formally, the functions can be defined as
follows:
%
%e3.51 #&#
\begin{equation}
\label{349} \qquad i(s) = \bigl(f(\cdot,s) - g(\cdot,s) \bigr)^{-1}(0)
\quad\mbox{and}\quad s(i) = \bigl(f(i,\cdot) - g(i,\cdot) \bigr)^{-1}(0)
\end{equation}
for $f(i,s) > g(i,s)$. [Recall from (\ref{334}) and
(\ref{337}) that $f(-1+,s) = -1$ and $g(-1+,s) < 1$ as well as that
$f(i,1-) > -1$ and $g(i,1-)=1$ for $-1 < i < s < 1$.] Geometrically,
moving from $i$ down to $i(s)$ (with $s$ fixed) corresponds to
moving along the first coordinate from any $(i,x,s)$ in $C_{f,g}^0$
to the closest point at the boundary between $C_{f,g}^0$ and
$C_{f,g}^-$ if $x \le f(i(s),s)$ and to\vspace*{1pt} the closest point at the
boundary between $C_{f,g}^0$ and $C_{f,g}^+$ if $x \ge f(i(s),s)$.
Similarly, moving from $s$ up to $s(i)$ (with $i$ fixed) corresponds
to moving along the third coordinate from any $(i,x,s)$ in
$C_{f,g}^0$ to the closest point at the\vspace*{1pt} boundary between $C_{f,g}^0$
and $C_{f,g}^+$ if $x \ge g(i,s(i))$ and to the closest point at the
boundary between $C_{f,g}^0$ and $C_{f,g}^-$ if $x \le g(i,s(i))$.

Since $(i(s),x,s)$ with $x \le f(i(s),s)$ belongs to the boundary of
$C_{f,g}^-$, we know that $V_{f,g}(i(s),x,s)$ is given by
(\ref{325}) above. Writing the integral from $x$ to $ f(i(s),s)$ in
this expression as the integral from $0$ to $f(i(s),s)$ minus the
integral from $0$ to $x$, it is easily seen that (\ref{325}) reads
as follows:
%
%e3.52 #&#
\begin{eqnarray}
\label{350}
V \bigl(i(s),x,s \bigr) &=& s - i(s)\nonumber
\\
&&{} + \bigl[c_2(s) -
c_1 \bigl(i(s) \bigr) \bigr]
\biggl[ H(x) - L(x) \int_0^{f(i(s),s)}
m(dy)
\\
&&\hspace*{123pt}{} + \int_0^{f(i(s),s)} L(y) m(dy) \biggr]\nonumber
\end{eqnarray}
for $x \le f(i(s),s)$. Comparing (\ref{350}) with (\ref{345}), we
can conclude that
%
%e3.53 #&#
%e3.54 #&#
\begin{eqnarray}
\label{351} \qquad a_1 \bigl(i(s) \bigr) + a_2(s) &=& -
\bigl[c_2(s) - c_1 \bigl(i(s) \bigr) \bigr] \int
_0^{
f(i(s),s)} m(dy),
\\
\label{352} b_1 \bigl(i(s) \bigr) + b_2(s) &=& s -
i(s) + \bigl[c_2(s) - c_1 \bigl(i(s) \bigr) \bigr] \int
_0^{f(i(s),s)} L(y) m(dy).
\end{eqnarray}
Using (\ref{346})--(\ref{347}) and (\ref{351})--(\ref{352}) we can
calculate $a_2'(s)$. First, by (\ref{347}) we can express
$a_2'(s)$ in terms of $b_2'(s)$. Second, by (\ref{352}) we can
express $b_2'(s)$ in terms of $b_1'(i(s))$. Third, by (\ref{346})
we can express $b_1'(i(s))$ in terms of $a_1'(i(s))$. Fourth, by
(\ref{351}) we can express $a_1'(i(s))$ in terms of $a_2'(s)$. This
closes the loop and gives an equation for $a_2'(s)$. A lengthy
calculation following these steps and making use of (\ref{36})
above yields
%
%e3.55 #&#
\begin{eqnarray}
\label{353} \qquad a_2'(s) &=& - \frac{1}{L(s) - L(i(s))}\nonumber
\\
&&{}\times  \biggl[
\frac{f_s'(i(s)
,s)}{f_i'(i(s),s)} \biggl[ 1 + c_1' \bigl(i(s) \bigr)
\int_{i(s)}^{f(i(s),s)} \bigl[L(y) - L \bigl(i(s) \bigr)
\bigr] m(dy) \biggr]
\\
&&\hspace*{32pt}{} + 1 + c_2'(s) \biggl[ H(s) + \int
_0^{f(i(s),s)} \bigl[L(y) - L \bigl(i(s) \bigr) \bigr]
m(dy) \biggr] \biggr].\nonumber
\end{eqnarray}
Similarly, since $(i,x,s(i))$ with $x \ge g(i,s(i))$ belongs to the
boundary of $C_{f,g}^+$ we know that $V_{f,g}(i,x,s(i))$ is given by
(\ref{332}) above. Writing the integral from $g(i,s(i))$ to $x$ in
this expression as the integral from~$0$ to $x$ minus the integral
from $0$ to $g(i,s(i))$, it is easily seen that (\ref{332}) reads
as follows:
%
%e3.56 #&#
\begin{eqnarray}
\label{354} V \bigl(i,x,s(i) \bigr) &=& s(i) - i \nonumber
\\
&&{} + \bigl[c_2
\bigl(s(i) \bigr) - c_1(i) \bigr]
\\
&&\quad {} \times \biggl[ H(x) - L(x) \int_0^{g(i,s(i))}
m(dy) + \int_0^{g(i,s(i))} L(y) m(dy) \biggr]\hspace*{-15pt}\nonumber
\end{eqnarray}
for $x \ge g(i,s(i))$. Comparing (\ref{354}) with (\ref{345}) we
can conclude that
%
%e3.57 #&#
%e3.58 #&#
\begin{eqnarray}
\label{355}\qquad a_1(i) + a_2 \bigl(s(i) \bigr) &=& -
\bigl[c_2 \bigl(s(i) \bigr) - c_1(i) \bigr] \int
_0^{
g(i,s(i))} m(dy),
\\
\label{356} b_1(i) + b_2 \bigl(s(i) \bigr) &=& s(i) -
i + \bigl[c_2 \bigl(s(i) \bigr) - c_1(i) \bigr] \int
_0^{g(i,s(i))} L(y) m(dy).
\end{eqnarray}
Using (\ref{346})--(\ref{347}) and (\ref{355})--(\ref{356}) we can
calculate $a_1'(i)$. First, by (\ref{346}) we can express
$a_1'(i)$ in terms of $b_1'(i)$. Second, by (\ref{356}) we can
express $b_1'(i)$ in terms of $b_1'(s(i))$. Third, by (\ref{347})
we can express $b_2'(s(i))$ in terms of $a_2'(s(i))$. Fourth, by
(\ref{355}) we can express $a_2'(s(i))$ in terms of $a_1'(i)$. This
closes the loop and gives an equation for $a_1'(i)$. A lengthy
calculation following these steps and making use of (\ref{37})
above yields
%
%e3.59 #&#
\begin{eqnarray}
\label{357} \qquad a_1'(i) &=& - \frac{1}{L(s(i)) - L(i)}\nonumber
\\
&&{}\hphantom{-}\times  \biggl[\frac{g_i'(i,
s(i))}{g_s'(i,s(i))} \biggl[ 1 + c_2' \bigl(s(i) \bigr)
\int_{g(i,s(i))}^{s(i)} \bigl[L \bigl(s(i) \bigr) - L(y)
\bigr] m(dy) \biggr]
\\
&&\hspace*{26.5pt}\hphantom{-}{} + 1 + c_1'(i) \biggl[ H(i) - \int
_0^{g(i,s(i))} \bigl[L \bigl(s(i) \bigr) - L(y) \bigr]
m(dy) \biggr] \biggr].\nonumber
\end{eqnarray}

We can now determine $A$ and $B$ in (\ref{341}) using the
closed-form expressions obtained. First, note that by (\ref{351})
we find that
%
%e3.60 #&#
\begin{eqnarray}
\label{358} A(i,s) &=& A \bigl(i(s),s \bigr) + \int_{i(s)}^i
A_u'(u,s) \,du\nonumber
\\
&=&  a_1 \bigl(i(s) \bigr) + a_2(s) + \int_{i(s)}^i a_1'(u) \,du
\nonumber\\[-8pt]\\[-8pt]
&=& - \bigl[c_2(s) - c_1 \bigl(i(s) \bigr)
\bigr] \int_0^{f(i(s),s)} m(dy)\nonumber
\\
&&{}  + \int _{i(s)}^i a_1'(u) \,du,\nonumber
\end{eqnarray}
where $a_1'(u)$ is given by (\ref{357}) above. Note also that by
(\ref{355}) we find that
%
%e3.61 #&#
\begin{eqnarray}
\label{359} A(i,s) &=& A \bigl(i,s(i) \bigr) - \int_s^{s(i)}
A_v'(i,v) \,dv\nonumber
\\
& =& a_1(i) + a_2
\bigl(s(i) \bigr) - \int_s^{s(i)}
a_2'(v) \,dv
\\
&=& - \bigl[c_2 \bigl(s(i) \bigr) - c_1(i)
\bigr] \int_0^{g(i,s(i))} m(dy) - \int
_s^{s(i)} a_2'(v) \,dv,\nonumber
\end{eqnarray}
where $a_2'(v)$ is given by (\ref{353}) above. Second, observe
that (\ref{346}) and (\ref{347}) yield
%
%e3.62 #&#
%e3.63 #&#
\begin{eqnarray}
\label{360} b_1'(i) &=& -a_1'(i)
L(i) + c_1'(i) H(i),
\\
\label{361} b_2'(s) &=& -a_2'(s)
L(s) - c_2'(s) H(s),
\end{eqnarray}
where $a_1'(i)$ and $a_2'(s)$ are given by (\ref{357}) and
(\ref{353}) above. Note that by (\ref{352}) we find that
%
%e3.64 #&#
\begin{eqnarray}
\label{362} B(i,s) &=& B \bigl(i(s),s \bigr) + \int_{i(s)}^i
B_u'(u,s) \,du\nonumber
\\
& =& b_1 \bigl(i(s) \bigr) +
b_2(s) + \int_{i(s)}^i
b_1'(u) \,du
\\
&=& s - i(s) + \bigl[c_2(s) - c_1 \bigl(i(s)
\bigr) \bigr] \int_0^{f(i(s),s)} L(y) m(dy) + \int _{i(s)}^i b_1'(u) \,du,\nonumber
\end{eqnarray}
where $b_1'(u)$ is given by (\ref{360}) above. Note also that by
(\ref{356}) we find that
%
%e3.65 #&#
\begin{eqnarray}
\label{363} B(i,s) &=& B \bigl(i,s(i) \bigr) - \int_s^{s(i)}
B_v'(i,v) \,dv\nonumber
\\
& =& b_1(i) + b_2
\bigl(s(i) \bigr) - \int_s^{s(i)}
b_2'(v) \,dv
\\
& =& s - i(s) + \bigl[c_2 \bigl(s(i) \bigr) -
c_1(i) \bigr] \int_0^{g(i,s(i))} L(y) m(dy) - \int_s^{s(i)} b_2'(v)\,dv,\nonumber
\end{eqnarray}
where $b_2'(v)$ is given by (\ref{361}) above.\vadjust{\goodbreak}

Finally, inserting (\ref{358}), (\ref{362}) and
(\ref{359}), (\ref{363}) into (\ref{341}) we, respectively, obtain
the following two closed-form expressions:
%
%e3.66 #&#
%e3.67 #&#
\begin{eqnarray}
\qquad V(i,x,s) &=& s - i(s) + \bigl[c_2(s) -
c_1 \bigl(i(s) \bigr) \bigr] \int_0^{f(i(s),s)}
\bigl[L(y) - L(x) \bigr] m(dy)\nonumber
\\
\label{364}  &&{} + \bigl[ c_2(s) - c_1(i) \bigr] H(x)
\\[-1pt]
&&{}  + \int _{i(s)}^i \bigl( \bigl[L(x) - L(u) \bigr]
a_1'(u) + c_1'(u) H(u)
\bigr) \,du,\nonumber
\\
V(i,x,s) &=& s(i) - i + \bigl[c_2 \bigl(s(i) \bigr) -
c_1(i) \bigr] \int_0^{g(i,
s(i))} \bigl[L(y)
- L(x) \bigr] m(dy)\nonumber
\\
\label{365} &&{} + \bigl[c_2(s) - c_1(i) \bigr] H(x)
\\[-1pt]
&&{} + \int _s^{s(i)} \bigl( \bigl[L(v) - L(x) \bigr]
a_2'(v) + c_2'(v) H(v) \bigr)\,dv\nonumber
\end{eqnarray}
for $f(i,s) > g(i,s)$ where $a_1'(u)$ and $a_2'(v)$ are given by
(\ref{357}) and (\ref{353}) above. A~formal verification of
(\ref{364}) and (\ref{365}) can be easily done by It\^o's formula
once we derive the optimality in the next step; see Remark~\ref{rem2} below.
Observe that if $f(i,s) = g(i,s)$, then $i(s)=i$ and $s(i)=s$ so that
the second integral in both (\ref{364}) and (\ref{365}) is zero,
and these expressions reduce to (\ref{325}) and (\ref{332}),
respectively.
\end{longlist}
\begin{longlist}[(3)]
\item[(7)] \emph{Optimality of the minimal and maximal solution}. We will
begin by disclosing the superharmonic characterisation of the value
function in terms of the solutions to (\ref{36}) and (\ref{37})
staying strictly above/below the lower/upper diagonal, respectively.
For this, let $i \mapsto f(i,s)$ be any solution to (\ref{36})
satisfying $f(i,s) > i$ for all $i$ with $f(-1+,s) \in[-1,1)$, and
let $s \mapsto g(i,s)$ be any solution to (\ref{37}) satisfying
$g(i,s) < s$ for all $s$ with $g(i,1-) \in(-1,1]$. Consider the
function $V_{f,g}$ defined by (\ref{325}) and (\ref{332}) on
$C_{f,g}^- \cup C_{f,g}^+$, and\vspace*{2pt} set $V_{f,g}(i,x,s) = s - i$ on
$D_{f,g}$ which denotes the complement of $C_{f,g}$. Then the same
arguments as in (\ref{335}) and (\ref{338}) above show that $s
\mapsto f(i,s)$ and $i \mapsto g(i,s)$ are decreasing. This implies
that after starting in the set $C_{f,g}^- \cup C_{f,g}^+ \cup
D_{f,g}$, the process $(I,X,S)$ remains in the same set for the rest
of time (i.e., it never enters the set $C_{f,g}^0$). Fix any
point $(i,x,s)$ such that $f(i,s) \le g(i,s)$ with $i \le x \le s$.
Note that $(i,x,s)$ belongs to $C_{f,g}^- \cup C_{f,g}^+ \cup
D_{f,g}$, and consider the motion of $(I,X,S)$ under $\mathsf P_{i,x,s}$.
Recall that $V_{f,g}$ solves the free boundary problem
(\ref{312})--(\ref{318}) on $C_{f,g}^- \cup C_{f,g}^+$. Due to the
``triple-deck'' structure of $V_{f,g}$ we can apply the
change-of-variable formula with local time on surfaces \cite{Pe-3}
which in view of (\ref{317}) and (\ref{318}) (note that these
conditions can fail for the second derivatives) reduces to standard
It\^o's formula and gives
%
%e3.68 #&#
\begin{eqnarray}\label{366}
&&  V_{f,g}(I_t,X_t,S_t)\nonumber
\\
&&\qquad = V_{f,g}(i,x,s) + \int_0^t
\frac{\partial
V_{f,g}}{\partial i}(I_s,X_s,S_s)
\,dI_s + \int_0^t
\frac{\partial
V_{f,g}}{\partial x}(I_s,X_s,S_s)
\,dX_s\hspace*{-25pt}
\\
&&\quad\qquad {} + \int_0^t \frac{\partial V_{f,g}}{\partial s}(I_s,X_s,S_s)
\,dS_s + \frac{1}{2} \int_0^t
\frac{\partial^2 V_{f,g}}{\partial x^2}(I_s,X_s,S_s) \,d\langle X,X
\rangle _s\nonumber
\\
&&\qquad = V_{f,g}(i,x,s) + \int_0^t
\sigma(X_s) \frac{\partial V_{f,g}}{\partial x}(I_s,X_s,S_s)
\,dB_s\nonumber
\\
&&\quad\qquad{} + \int_0^t (\mathbb{L}_X
V_{f,g}) (I_s,X_s,S_s) \,ds,\nonumber
\end{eqnarray}
where we also use (\ref{313}) and (\ref{314}) to conclude that the
integrals with respect to $dI_s$ and $dS_s$ are equal to zero. The
process $M=(M_t)_{t \ge0}$ defined by
%
%e3.69 #&#
\begin{equation}
\label{367} M_t = \int_0^t
\sigma(X_s) \frac{\partial V_{f,g}}{\partial x} (I_s,X_s,S_s)
\,dB_s
\end{equation}
is a continuous local martingale. Introducing the increasing process
$P=(P_t)_{t \ge0}$ by setting
%
%e3.70 #&#
\begin{equation}
\label{368} P_t = \int_0^t
c(I_s,X_t,S_s) 1 \bigl(f(I_s,S_s)
\le X_s \le g(I_s,X_s) \bigr) \,ds
\end{equation}
and using the fact that the set of all $s$ for which $X_s$ is either
$f(I_s,S_s)$ or $g(I_s,S_s)$ is of Lebesgue measure zero, we see by
(\ref{312}) that (\ref{366}) can be rewritten as follows:
%
%e3.71 #&#
\begin{equation}
\label{369} V_{f,g}(I_t,X_t,S_t)
- \int_0^t c(I_s,X_s,S_s)
\,ds = V_{f,g} (i,x,s) + M_t - P_t.
\end{equation}
From this representation we see that the process
\[
V_{f,g}(I_t,X_t,S_t) - \int_0^t c(I_s,X_s,S_s) \,ds
\]
is a local
supermartingale for $t \ge0$.

Let $\tau$ be any stopping time of $X$. Choose a localisation
sequence $(\sigma_n)_{n \ge1}$ of bounded stopping times for $M$.
From (\ref{325}) and (\ref{332}) we see that $V_{f,g}(i,x,s) \ge s -
i$ for all $(i,x,s) \in C_{f,g}^- \cup C_{f,g}^+ \cup D_{f,g}$.
Recalling that the process $(I,X,S)$ remains in the latter set, we
can conclude from (\ref{369}) using the optional sampling theorem
that
%
%e3.72 #&#
\begin{eqnarray}
\label{370} &&\mathsf E _{i,x,s} \biggl[ S_{\tau\wedge\sigma_n} -
I_{\tau\wedge
\sigma_n} - \int_0^{\tau\wedge\sigma_n}
c(I_s,X_s,S_s) \,dt \biggr]\nonumber
\\
&&\qquad \le\mathsf E _{i,x,s} \biggl[ V_{f,g}(I_{\tau
\wedge\sigma_n},X_{\tau\wedge\sigma_n},S_{\tau\wedge\sigma_n})
- \int_0^{\tau\wedge\sigma_n} c(I_s,X_s,S_s)
\,dt \biggr]
\\
\nonumber
&&\qquad \le V_{f,g}(i,x,s) + \mathsf E _{i,x,s}(M_{\tau
\wedge\sigma_n})
= V_{f,g}(i,x,s)
\end{eqnarray}
for all $(i,x,s) \in C_{f,g}^- \cup C_{f,g}^+ \cup D_{f,g}$ and all
$n \ge1$. Letting $n \rightarrow\infty$ and using the monotone
convergence theorem we find that
%
%e3.73 #&#
\begin{equation}
\label{371} \mathsf E _{i,x,s} \biggl[ S_\tau -
I_\tau- \int_0^\tau
c(I_s,X_s, S_s) \,dt \biggr] \le
V_{f,g}(i,x,s)
\end{equation}
for all $(i,x,s) \in C_{f,g}^- \cup C_{f,g}^+ \cup D_{f,g}$. Taking
first the supremum over all $\tau$ and then the infimum over all
$f$ and $g$, we conclude that
%
%e3.74 #&#
\begin{equation}
\label{372} V(i,x,s) \le\inf_{f,g} V_{f,g}(i,x,s)
= V_{f_*,g_*}(i,x,s)
\end{equation}
for all $(i,x,s) \in C_{f,g}^- \cup C_{f,g}^+ \cup D_{f,g}$ where
$f_*$ denotes the minimal solution to (\ref{36}) staying strictly
above the lower diagonal, and $g_*$ denotes the maximal solution to~(\ref{37}) staying strictly below the upper diagonal. Recalling
that $f \mapsto V_{f,g}$ is increasing and $g \mapsto V_{f,g}$ is
decreasing when $f \le g$, we see that the infimum in (\ref{372}) is
attained over any sequence of solutions $f_n$ and $g_n$ to
(\ref{36}) and (\ref{37}) such that $f_n \downarrow f_*$ and $g_n
\uparrow g_*$ as $n \rightarrow\infty$. Since $f_*$ and $g_*$ are
solutions themselves to which (\ref{371}) applies, we see that
(\ref{372}) holds for all $(i,x,s)$ in the set $C_{f_*,g_*}^- \cup
C_{f_*,g_*}^+ \cup D_{f_*,g_*}$ which is the increasing union of the
sets $C_{f_n,g_n}^- \cup C_{f_n,g_n}^+ \cup D_{f_n,g_n}$ for $n \ge
1$. From these considerations and (\ref{372}) in particular, it
follows that the only possible candidates for the optimal stopping
boundary are the minimal and maximal solution $f_*$ and $g_*$. Note
that (\ref{370}) also implies that
%
%e3.75 #&#
\begin{equation}
\label{373} \mathsf E _{i,x,s} \biggl[ V_{f,g}(I_\tau,X_\tau,S_\tau)
- \int_0^\tau c(I_s,X_s,S_s)
\,dt \biggr] \le V_{f,g}(i,x,s)
\end{equation}
showing that the function $(i,x,s,a) \mapsto V_{f,g}(i,x,s) - a$ is
superharmonic for the Markov process $(I,X,S,A)$ on the set
$C_{f,g}^- \cup C_{f,g}^+ \cup D_{f,g}$ where $A_t = \int_0^t
c(I_s,X_s,S_s) \,ds$ for $t \ge0$. Recalling that $f \mapsto
V_{f,g}$ is increasing and $g \mapsto V_{f,g}$ is decreasing when $f
\le g$, and that $V_{f,g}(i,x,s) \ge s - i$ for all $(i,x,s) \in
C_{f,g}^- \cup C_{f,g}^+ \cup D_{f,g}$, we see\vspace*{1pt} that selecting the
minimal solution $f_*$ staying strictly above the lower diagonal and
the maximal solution $g_*$ staying strictly below the upper diagonal
is equivalent to invoking the superharmonic characterisation of the
value function (according to which the value function is the
smallest superharmonic function which dominates the gain function).
For more details on the latter characterisation in a general setting
we refer to \cite{PS}, Chapter~1; see also Remark~\ref{rem3} below.

To\vspace*{1pt} prove that $f_*$ and $g_*$ are optimal on $ C_{f_*,g_*}^- \cup
C_{f_*,g_*}^+ \cup D_{f_*,g_*}$, consider the stopping time
$\tau_{f_n,g_n}$ defined in (\ref{38}) where $i \mapsto f_n(i,s)$
is the solution to~(\ref{36}) such that $f_n(i_n,s) = i_n$ and $s
\mapsto g_n(i,s)$ is the solution to (\ref{37}) such that
$g_n(i,s_n) = s_n$ for some $i_n \downarrow-1$ and $s_n \uparrow1$
as $n \rightarrow\infty$. Consider the function $V_{f_n,g_n}$
defined by (\ref{325}) and (\ref{332}) on $C_{f_n,g_n}^- \cup
C_{f_n,g_n}^+$, and\vspace*{1pt} set $V_{f_n,g_n}(i,x,s) = s - i$ for $(i,x,s)
\in D_{f_n,g_n}$ and $n \ge1$. Recall that $V_{f_n,g_n}$ solves the
free-boundary problem \mbox{(\ref{312})--(\ref{318})} on $C_{f_n,g_n}^-
\cup C_{f_n,g_n}^+$ for $n \ge1$. Fix any $(i,x,s)$ in $
C_{f_*,g_*}^- \cup C_{f_*,g_*}^+ \cup D_{f_*,g_*}$, and note that
this $(i,x,s)$ belongs to $ C_{f_n,g_n}^- \cup C_{f_n,g_n}^+ \cup
D_{f_n,g_n}$ since $f_n \le f_*$ and $g_* \le g_n$ for every $n \ge
1$. The same arguments as above yield the formula~(\ref{366}) with
$f_n$ and $g_n$ in place of $f$ and $g$ for $n \ge1$. Since
$\sigma$ and $\partial V_{f_n,g_n}/\partial x$ are bounded on
$C_{f_n,g_n}^- \cup C_{f_n,g_n}^+$, we see\vspace*{1pt} that $(M_{t \wedge
\tau_{f_n,g_n}})_{t \ge0}$ defined by (\ref{367}) with $f_n$ and
$g_n$ in place of $f$ and $g$ is a martingale under $\mathsf P_{ i,x,s}$.
The latter conclusion follows from the fact that $\tau_{f_n,g_n} \le
\rho_{i_n,s_n}$ with $\mathsf E _{i,x,s} \rho_{i_n,s_n} < \infty$
implying also that $\mathsf E _{i,x,s} \int_0^{\tau_{f_n,g_n}}
c(I_s,X_s,S_s) \,dt < \infty$ for $n \ge1$. Since the process $P$
defined by (\ref{368}) with $f_n$ and $g_n$ in place of $f$ and $g$
satisfies $P_{\tau_{f_n,g_n}} = 0$, it follows from (\ref{369})
using (\ref{315}) and (\ref{316}) that
%
%e3.76 #&#
\begin{equation}
\label{374} \qquad V_{f_n,g_n}(i,x,s) = \mathsf E _{i,x,s} \biggl[
S_{\tau_{f_n,g_n}} - I_{\tau_{f_n,g_n}} - \int_0^{\tau_{f_n,g_n}}
c(I_s,X_s, S_s) \,dt \biggr]
\end{equation}
for all $i \le x \le s$ such that $f_n(i,s) \le g_n(i,s)$ with $n
\ge1$. Letting $n \rightarrow\infty$ in (\ref{374}), noting that
$\tau_{f_n,g_n} \uparrow\tau_{f_*,g_*}$ (since $[-1,1]^3$ is
compact), and using the monotone convergence theorem (recalling that
$S_{\tau_{f_*,g_*}} - I_{\tau_{f_*,g_*}}$ is bounded by $2$ and
therefore integrable) we find that
%
%e3.77 #&#
\begin{equation}
\label{375} \qquad V_{f_*,g_*}(i,x,s) = \mathsf E _{i,x,s} \biggl[
S_{\tau_{f_*,g_*}} - I_{\tau_{f_*,g_*}} - \int_0^{\tau_{f_*,g_*}}
c(I_s,X_s, S_s) \,dt \biggr]
\end{equation}
for all $i \le x \le s$ such that $f_*(i,s) \le g_*(i,s)$. This
shows that we have equality in (\ref{372}) and completes the proof
of the optimality of $\tau_{f_*,g_*}$ on the set $ C_{f_*,g_*}^-
\cup C_{f_*,g_*}^+ \cup D_{f_*,g_*}$.

To prove the optimality of $\tau_{f_*,g_*}$ on the set $
C_{f_*,g_*}^0$, that is, when $f_*(i,s) > g_*(i,s)$ for some $(i,s)$
given and fixed, one could attempt to apply similar arguments to
those in (\ref{370}) above. For this, however, we would need to
know that $V_{f,g}(i,x,s) \ge s - i$ not only for $f(i,s) \le
g(i,s)$ as follows from the closed-form expressions
(\ref{325}) and (\ref{332}) above but also for $f(i,s) > g(i,s)$. A
closer inspection of the latter case indicates that this
verification may be problematic if it is to follow from similar
closed-form expressions. Indeed, even in the special case of $c(i,s)
= c_2(s) - c_1(i)$, we see from (\ref{364}) and (\ref{365}) that
the conclusion is unclear since $a_1'(u)$ and $a_2'(v)$ appearing
there could also (at least in principle) take negative values as
well; see (\ref{353}) and (\ref{357}) above. To overcome this
difficulty we will exploit the extremal properties of the candidate
surfaces $f_*$ and $g_*$ in an essential way (in many ways this can
be seen as a key argument in the proof showing the full power of the
method). For this, take any point $(i_0,x_0,s_0)$ in the state space
such that $f_*(i_0,s_0) > g_*(i_0,s_0)$ with $i_0 < x_0 < s_0$ and
fix any $d_0 \in(i_0 \vee g_*(i_0,s_0),s_0 \wedge f_*(i_0,s_0))
\setminus\{x_0\}$. Choose solutions $i \mapsto f_d(i,s_0)$ and $s
\mapsto g_d(i_0,s)$ to (\ref{36}) and (\ref{37}) such that
$f_d(i_0,s_0) = d_0$ and $g_d(i_0,s_0) = d_0$, respectively. Note
that this is possible since $d_0$ lies strictly between $i_0$ and
$f_*(i_0,s_0)$ in the first case and strictly between $g_*(i_0,s_0)$
and $s_0$ in the second case. Note also that $i \mapsto f_d(i,s_0)$
must hit the lower diagonal and $s \mapsto g_d(i_0,s)$ must hit the
upper diagonal since $i \mapsto f_*(i,s_0)$ and $s \mapsto
g_*(i_0,s)$ are the minimal and maximal solutions staying strictly
above/below the lower/upper diagonal, respectively. Moreover, by the
construction of $f_d$ and $g_d$ we see that $(i_0,x_0,s_0)$ belongs
to either $C_{f_d,g_d}^-$ if $x_0 < d_0$ or $C_{f_d,g_d}^+$ if $x_0
> d_0$, and after starting at $(i_0,x_0,s_0)$ the process $(I,X,S)$
remains in either $C_{f_d,g_d}^-$ or $C_{f_d,g_d}^+$, respectively,
before hitting $D_{f_d,g_d}$. Considering the stopping time
$\tau_{f_d,g_d}$ defined in (\ref{38}) we therefore see that the
same arguments as those leading to (\ref{374}) also show that
%
%e3.78 #&#
\begin{eqnarray}\label{376}
&& V_{f_d,g_d}(i_0,x_0,s_0)
\nonumber\\[-8pt]\\[-8pt]
&&\qquad = \mathsf E _{i_0,x_0,s_0} \biggl[ S_{\tau_{f_d,
g_d}} - I_{\tau_{f_d,g_d}} - \int
_0^{\tau_{f_d,g_d}} c(I_s,X_s,S_s)
\,dt \biggr],\nonumber
\end{eqnarray}
where $V_{f_d,g_d}$ is given by either (\ref{325}) or (\ref{332}),
respectively. From the latter closed-form expressions we see that
$V_{f_d,g_d}(i_0,x_0,s_0) > s_0 - i_0$ and from (\ref{376}) it
therefore follows that $(i_0,x_0,s_0)$ belongs to the continuation
set $C$. Combining this conclusion\vspace*{-1pt} with the description of the
stopping set $D$ outside $C_{f_*,g_*}^0$ derived above, we see that
$C = C_{f_*,g_*}^0 \cup C_{f_*,g_*}^- \cup C_{f_*,g_*}^+$. This\vspace*{2pt}
proves the optimality of $\tau_*$ in (\ref{35}) and completes the
proof.\quad\qed
\end{longlist}\noqed
\end{pf}

We conclude this section with a few remarks on the preceding result
and proof.

%re1 #&#
\begin{rem}\label{rem1}
To describe the nature of the optimal stopping
time $\tau_{f_*,g_*}$ from~(\ref{35}), assume that the process
$(I,X,S)$ starts at $(0,0,0)$. Then due to $g_*(0,0) < 0 < f_*(0,0)$
we see that it is not optimal to stop at once so that $t \mapsto
I_t$ and $t \mapsto S_t$ will gradually start to decrease and
increase whenever $t \mapsto X_t$ returns to the lower and upper
diagonal, respectively. Due to (\ref{334})--(\ref{335}) and
\mbox{(\ref{337})--(\ref{338})} we see that $t \mapsto f_*(I_t,S_t)$ is
decreasing and $t \mapsto g_*(I_t,S_t)$ is increasing. Since $I_t
\downarrow-1$ and/or $S_t \uparrow1$ as $t \uparrow\infty$ we see
from (\ref{334}) and (\ref{337}) that the two sample paths $t
\mapsto f_*(I_t,S_t)$ and $t \mapsto g_*(I_t,S_t)$ will meet at some
random time which coincides with the first exit time of $(I,X,S)$
from\vspace*{1pt} the set $C_{f_*,g_*}^0$ defined in (\ref{39}). This can only
happen either through the lower diagonal (when $X$ is equal to $I$)
or through the upper diagonal (when $X$ is equal to $S$). In the
former case the process $(I,X,S)$ enters the set $C_{f_*,g_*}^-$
defined in (\ref{310}) and in the latter case the process $(I,X,S)$
enters the set $C_{f_*,g_*}^+$ defined in (\ref{311}). After
entering either $C_{f_*,g_*}^-$ or $C_{f_*,g_*}^+$ the process
$(I,X,S)$ remains in the same set until the first hitting of $X$ to
either $f(I,S)$ from below or $g(I,S)$ from above happens,
respectively. This moment defines the optimal stopping time
$\tau_{f_*,g_*}$. Note that from the optimality derived in the proof
of Theorem~\ref{teo1} [recall (\ref{374}) and (\ref{375}) in particular] we see
that $\tau_{f_*,g_*}$ has finite expectation (since otherwise the
value function would be equal to $-\infty$ and as such
$\tau_{f_*,g_*}$ could not be optimal). Note that the analogous
description of $\tau_{f_*,g_*}$ also holds for any starting point
$(i,x,s)$ of $(I,X,S)$ in the state space. After starting in
$C_{f_*,g_*}^0$ the process $(I,X,S)$ enters either $C_{f_*,g_*}^-$
or $C_{f_*,g_*}^+$ to remain in the same set until $\tau_{f_*,g_*}$
happens. The latter fact also holds if $(I,X,S)$ starts in either
$C_{f_*,g_*}^-$ or $C_{f_*,g_*}^+$ directly. To visualise these
movements, see Figure~\ref{fig1} above and note that $i_0$ and $s_0$ mark the
borderline levels between $C_{f_*,g_*}^0$ and $C_{f_*,g_*}^- \cup
C_{f_*,g_*}^+$ as described above.
\end{rem}

%re2 #&#
\begin{rem}\label{rem2}
Although we do not make use of this fact in the
proof of the optimality above, we note that in addition to the
closed-form expressions (\ref{325}) and (\ref{332}) on $C_{f,g}^-$
and $C_{f,g}^+$, respectively, the probabilistic representation
(\ref{319}) itself can also be used to define the function
$V_{f,g}$ on $C_{f,g}^0$ when the stopping time $\tau_{f,g}$ from
(\ref{38}) has finite expectation, and the resulting function will
solve the free boundary problem (\ref{312})--(\ref{318}) on
$C_{f,g}$ for the surfaces $f$ and $g$ constructed in the proof
above (those hitting the lower/upper diagonal at a single point and
the minimal/maximal solutions staying above/below the lower/upper
diagonal). Indeed, due to the monotonicity properties of $f$ and $g$
derived above, we see that after starting in $C_{f,g}^0$, the process
$(I,X,S)$ enters either the set $C_{f,g}^-$ or the set $C_{f,g}^+$
through the boundary $f=g$ to stay in the same set until
$\tau_{f,g}$ happens. This shows that defining the function
$V_{f,g}$ by (\ref{319}) on $C_{f,g}^0$ corresponds to solving the
Dirichlet problem stochastically where the value at the boundary
$f=g$ is set to be either (\ref{325}) at the lower diagonal or
(\ref{332}) at the upper diagonal, respectively. For standard
arguments how this can be done including how the required smoothness
of $V_{f,g}$ on $C_{f,g}^0$ can be derived; see, for example, \cite{PS}, Sections~7.1--7.3.
\end{rem}

%re3 #&#
\begin{rem}\label{rem3}
In addition to the facts used in the proof above
it is also useful to know that the superharmonic characterisation of
the value function represents the ``dual problem'' to the primal
problem (\ref{34}). For more details on the meaning of this claim
including connections to the Legendre transform, see \cite{Pe-4}.
\end{rem}

%re4 #&#
\begin{rem}\label{rem4}
A closer look into the proof above indicates that
the arguments developed and/or used should be applicable in more
general settings of the optimal stopping problem (\ref{34}) and its
relatives. As stated above it is not essential that the state space
of the diffusion process $X$ equals $(-1,1)$, and the result and
methodology of Theorem~\ref{teo1} should be valid for more general state
spaces (including $\mathbb{R}$ and $\mathbb{R}_+$ in particular). In
this case we
may need to take the supremum in (\ref{34}) over all stopping times
such that the expectation of the integral is finite, and although
the stopping time $\tau_{f_*,g_*}$ may not belong to this class in
some particular examples [so that the right-hand side of (\ref{34})
may not even be well-defined], this stopping time should be
approximately optimal in the sense that the approximate stopping
times $\tau_{f_n,g_n}$ yield the value (\ref{34}) in the limit as
$n \rightarrow\infty$. These extensions also include various
boundary behaviour of the process $X$ at the endpoints of the state
space (e.g., $0$ when the state space equals $\mathbb{R}_+$). We leave
precise formulations of these statements and proofs as informal
conjectures open for future developments. We emphasise that these
questions are best studied through examples, and each particular
example may have specifics which are difficult to cover by any
meta-theorem in advance. Omitting further details we briefly turn to
some examples.
\end{rem}

%s4 #&#
\section{Examples}\label{sec4}
%%%%%%%%%%%%%%%%%%

Combining the results of Proposition~\ref{prop1} and Theorem~\ref{teo1} we obtain the
solution to the quickest detection problem (\ref{23}). We
illustrate various special cases of this correspondence through one
particular example.

%ex1 #&#
\begin{exa}\label{exa1}
Assume that the observed process $Z$ is a
standard Brownian motion $B$ starting at $0$, suppose that $\ell$ is
a standard normal random variable independent from $B$, and consider
the quickest detection problem (\ref{23}) where $c>0$ is a given
and fixed constant. By the result of Proposition~\ref{prop1} we know that this
problem is equivalent to the optimal stopping problem (\ref{214})
where $X = 2F(Z) - 1$ solves~(\ref{210}) with $\mu$ and $\sigma$
given by (\ref{211}) and (\ref{212}). From (\ref{21}) we see that
$a=0$ and $b=1$ so that
%
%e4.1 #&#
%e4.2 #&#
\begin{eqnarray}
\label{41} \mu(x) &=& -\Phi^{-1} \biggl( \frac{x+1}{2}
\biggr) \varphi \biggl( \Phi^{-1} \biggl( \frac{x+1}{2} \biggr)
\biggr),
\\
\label{42} \sigma(x) &=& 2 \varphi \biggl( \Phi^{-1} \biggl(
\frac{x+1}{2} \biggr) \biggr)
\end{eqnarray}
for $x \in(-1,1)$ where $\Phi(y) = (1/\sqrt{2 \pi})
\int_{-\infty}^y e^{-z^2/2} \,dz$ is the standard normal
distribution function and $\varphi(y) = (1/\sqrt{2 \pi})
e^{-y^2/2}$ is the standard normal density function for $y \in\mathbb{R}$.
It is easily verified using (\ref{216}) that the scale function of
$X$ can be taken as
%
%e4.3 #&#
\begin{equation}
\label{43} L(x) = \int_0^x \exp \biggl(
\frac{1}{2} \biggl( \Phi^{-1} \biggl( \frac{y+1}{2} \biggr)
\biggr)^{ 2} \biggr) \,dy
\end{equation}
for $x \in(-1,1)$. By Theorem~\ref{teo1} we know that the following stopping
time is optimal:
%
%e4.4 #&#
\begin{equation}
\label{44} \tau_* = \inf \bigl\{ t \ge0 \vert f_*(I_t,S_t)
\le X_t \le g_*(I_t,S_t) \bigr\},
\end{equation}
where the surfaces $f_*$ and $g_*$ are the minimal and maximal
solutions to
%
%%%{\fontsize{10.5pt}{11pt}\selectfont{\begin{eqnarray}
%%%&& \frac{\partial f}{\partial i}(i,s)\nonumber\hspace*{-18pt}
%%%\\
%%%\label{45}  &&\qquad = \frac{ 2 \varphi^2
%%%( \Phi^{-1}  ((f(i,s) + 1)/2  )  ) \exp
%%%( 1/2  (\Phi^{-1}  ((f(i,s) + 1)/2
%%%)  )^{ 2}  ) }{ c(s - i) \int_i^{f(i,s)} \exp
%%%( 1/2  (\Phi^{-1}  ((y + 1)/2  ))^{ 2}  ) \,dy }\hspace*{-18pt}
%%%\\
%%%&&\quad\qquad{} \times \biggl[ 1 + c \int_i^{f(i,s)}
%%%\frac{\int_i^y \exp ( 1/2
%%%(\Phi^{-1}  ( (z + 1)/2  )  )^{ 2}  )
%%%\,dz}{ 2 \varphi^2  ( \Phi^{-1}  ((y + 1)/2  )
%%%) \exp ( 1/2  (\Phi^{-1}  ((y +1)/2  )  )^{ 2}  ) } \,dy \biggr],\nonumber\hspace*{-18pt}
%%%\\
%%%&& \frac{\partial g}{\partial s}(i,s)\nonumber\hspace*{-18pt}
%%%\\
%%%\label{46}  &&\qquad = \frac{ 2 \varphi^2  (
%%%\Phi^{-1}  ((g(i,s) + 1)/2  )  ) \exp (
%%%1/2  (\Phi^{-1}  ( (g(i,s) + 1)/2  )
%%%)^{ 2}  ) }{ c(s - i) \int_{g(i,s)}^s \exp (
%%%1/2  (\Phi^{-1}  ((y + 1)/2 ))^{ 2}  ) \,dy }\hspace*{-18pt}
%%%\\
%%%&&\quad\qquad{} \times \biggl[ 1 + c \int_{g(i,s)}^s
%%%\frac{\int_y^s \exp ( 1/2
%%%(\Phi^{-1}  ( (z + 1)/2  )  )^{ 2}  )
%%%\,dz}{ 2 \varphi^2  ( \Phi^{-1}  ((y + 1)/2  )
%%%) \exp ( 1/2  (\Phi^{-1}  ((y +1)/2  )  )^{ 2}  ) } \,dy \biggr]\nonumber\hspace*{-18pt}
%%%\end{eqnarray}}}%
%e4.5 #&#
%e4.6 #&#
\begin{eqnarray}
&& \frac{\partial f}{\partial i}(i,s)\nonumber\hspace*{-6pt}
\\
\label{45}  &&\!\qquad = \frac{ 2 \varphi^2
( \Phi^{-1}  ((f(i,s) + 1)/2  )  ) \exp
( 1/2  (\Phi^{-1}  ((f(i,s) + 1)/2
)  )^{ 2}  ) }{ c(s - i) \int_i^{f(i,s)} \exp
( 1/2  (\Phi^{-1}  ((y + 1)/2  ))^{ 2}  ) \,dy }\hspace*{-6pt}
\\
&&\hspace*{28pt}{} \times \biggl[ 1 + c \int_i^{f(i,s)}\!
\frac{\int_i^y \exp ( 1/2
(\Phi^{-1}  ( (z + 1)/2  )  )^{ 2}  )
\,dz}{ 2 \varphi^2  ( \Phi^{-1}  ((y + 1)/2  )
) \exp ( 1/2  (\Phi^{-1}  ((y +1)/2  )  )^{ 2}  ) } \,dy \biggr],\nonumber\hspace*{-6pt}
\\
&& \frac{\partial g}{\partial s}(i,s)\nonumber\hspace*{-6pt}
\\
\label{46}  &&\!\qquad = \frac{ 2 \varphi^2  (
\Phi^{-1}  ((g(i,s) + 1)/2  )  ) \exp (
1/2  (\Phi^{-1}  ( (g(i,s) + 1)/2  )
)^{ 2}  ) }{ c(s - i) \int_{g(i,s)}^s \exp (
1/2  (\Phi^{-1}  ((y + 1)/2 ))^{ 2}  ) \,dy }\hspace*{-6pt}
\\
&&\hspace*{28pt} {} \times \biggl[ 1 + c \int_{g(i,s)}^s
\frac{\int_y^s \exp ( 1/2
(\Phi^{-1}  ( (z + 1)/2  )  )^{ 2}  )
\,dz}{ 2 \varphi^2  ( \Phi^{-1}  ((y + 1)/2  )
) \exp ( 1/2  (\Phi^{-1}  ((y +1)/2  )  )^{ 2}  ) } \,dy \biggr]\nonumber\hspace*{-6pt}
\end{eqnarray}
staying strictly above the lower diagonal $d^s$ and strictly below
the upper diagonal $d_i$ for $i<s$ in $(-1,1)$, respectively. Equations
(\ref{45}) and (\ref{46}) are singular at the lower and
upper diagonal. Passing to the inverse equations $\partial i /
\partial f$ and $\partial s /
\partial g$ these singularities get removed, and one can determine
the minimal and maximal solution by approximating them with the
solutions which hit the lower and upper diagonal, respectively (as
explained in the proof above). The results of these calculations are
illustrated in Figures~\ref{fig1}--\ref{fig3}. Similar qualitative behaviour of the
optimal surfaces can also be observed in other examples of
diffusions and hidden levels.

The list of examples can be continued by considering various
diffusion processes $Z$ and hidden targets $\ell$. This leads to a
classification of the laws of $\ell$ against the laws of $Z$
(through the drift and diffusion coefficient) in terms of the
optimal surfaces derived in Theorem~\ref{teo1}. This classification can be
used for calibration against observed performance (where either of
the two laws is taken initially to be known, e.g.).

Apart from the problems where the optimal stopping boundaries are
surfaces, this also includes problems where the optimal stopping
boundaries are curves. We illustrate this briefly through one-known
example from stochastic analysis.
\end{exa}

%ex2 #&#
\begin{exa}\label{exa2}
Taking $X$ to be a standard Brownian motion $B$
and setting $c(r) \equiv c$, it is easily seen that the minimal and
maximal solutions to (\ref{26})~and~(\ref{27}) are given by
%
%e4.7 #&#
\begin{equation}
\label{47} f(i,s) = i + \frac{1}{2c} \quad\mbox{and}\quad g(i,s) = s
- \frac
{1}{2c}.
\end{equation}
From (\ref{348}) we see that $i(s) = s - \frac{1}{c}$ and $s(i) =
i + \frac{1}{c}$. Since $f_s' \equiv0$ we see from (\ref{353})
that $a_2'(s) \equiv-c$. Inserting this into (\ref{365}) we find
that $V(0,0,0) = \frac{3}{4c}$; note that unboundedness of $B$
presents no difficulty since the optimal stopping time has finite
expectation. This shows that for any stopping time $\tau$ of $B$
(with finite expectation) we have
%
%e4.8 #&#
\begin{equation}
\label{48} \mathsf E (S_\tau - I_\tau) \le c \mathsf E
\tau+ \frac{3}{4c}.
\end{equation}
Taking the infimum over all $c>0$ we obtain the result of
\cite{DGM},
%
%e4.9 #&#
\begin{equation}
\label{49} \mathsf E (S_\tau - I_\tau) \le\sqrt{3}
\sqrt{\mathsf E \tau}.
\end{equation}
One can extract similar other inequalities/information from the
proof above.
\end{exa}


%suskaldyti doi

% imsref loaded by linak, 2014-03-03 12:29:08
% imsref loaded by linak, 2014-03-03 13:54:13
% imsref loaded by linak, 2014-03-04 14:32:06

\printaddresses

\end{document}